\documentclass[12pt]{article}
\usepackage[utf8]{inputenc}
\usepackage[margin=1in]{geometry}
\usepackage[sorting=none, style=numeric-comp]{biblatex}
\usepackage{amsmath}
\usepackage{amssymb}
\usepackage{setspace}
\usepackage{bm}
\usepackage{geometry}
\usepackage{dirtytalk}
\usepackage{graphicx}
\usepackage{algorithm}
\usepackage{algpseudocode}
\usepackage{xcolor}
\usepackage{hyperref}
\usepackage[capitalise]{cleveref}
\usepackage{multirow}
\usepackage[affil-it]{authblk}

\newtheorem{theorem}{Theorem}

\title{A Multidimensional Fourier Approximation of Optimal Control Surfaces}

\author{Gabriel Nicolosi\thanks{Electronic address: \texttt{gabrielnicolosi@mst.edu}; Corresponding author}}
\affil{Dept.\ of Engineering Management and Systems Engineering, Missouri University of Science and Technology\\ Rolla, MO.}

\author{Christopher Griffin}
\affil{Applied Research Laboratory, Penn State University\\ University Park, PA.}

\author{Terry Friesz}
\affil{Harold and Inge Marcus Dept.\ of Industrial and Manufacturing Engineering, Penn State University\\ University Park, PA.}

\bibliography{References} 

\begin{document}

\maketitle

\begin{abstract}
    This work considers the problem of approximating initial condition and time-dependent optimal control and trajectory surfaces using multivariable Fourier series. A modified Augmented Lagrangian algorithm for translating the optimal control problem into an unconstrained optimization one is proposed and two problems are solved: a quadratic control problem in the context of Newtonian mechanics, and a control problem arising from an odd-circulant game ruled by the replicator dynamics. Various computational results are presented. Use of automatic differentiation is explored to circumvent the elaborated gradient computation in the first-order optimization procedure. Furthermore, mean square error bounds are derived for the case of one and two-dimensional Fourier series approximations, inducing a general bound for problems with state space of $n$ dimensions. 
\end{abstract}

\section{Introduction}
Recently, the fields of nonlinear dynamics and nonlinear control have witnessed an increasing interest in the application of machine learning methods to tackle the common analytically challenging problems pertaining to them. These models are mostly useful when the so-called curse of dimensionality hinders the efficient derivation of closed-form or approximated solutions to problems presenting nonlinearities and high dimensionality. To name a few, see \cite{BEG19, KMAWJ21, H16, ROLNF20, CYJD18, B22, KSW21, PY21}. This work is more generally considered in the context of a relevant sub-area of nonlinear systems, namely the class of evolutionary games \cite{W97, HS98} and, more specifically, odd-circulant games \cite{G16}. In \cite{GF22}, an optimal control problem (OCP) arising from odd-circulant games, such as rock-paper-scissors, is studied. In this work, we address the problem of finding an approximating continuous-time control surface for an OCP arising from such odd-circulant evolutionary games. The main contributions of this work are: 
\begin{enumerate}
    \item Inspired by some of the state-of-art numerical approaches for solving OCPs, we propose a Fourier-based, initial condition and time-dependent optimal control and trajectory surface approximations, seeking to generalize the approach originally introduced in \cite{N22, NG23};
    \item Building upon a recent effort to look at these numerical methods through a machine learning lens and, more specifically, neural networks (see \cite{N17, Y23, EP13, S18, KK12, BZ20, N23}), we adapt an Augmented Lagrangian algorithm, borrowed from the theory of nonlinear optimization, capable of translating the OCP into an unconstrained optimization problem suitable to the application of differentiation algorithms enabled by modern automatic differentiation software;
    \item Extending from the results in \cite{G72}, we derive mean square error (MSE) upper bounds for the truncated Fourier approximations in the single and double variable cases, obtaining theoretical estimates for the number of coefficients in the truncated approximations that guarantee a desired MSE threshold is not violated. Our results suggest that a conjecture for the $n-$dimensional case can be drawn; 
    \item Finally, we show how the multidimensional truncated Fourier series, allied with the modified Augmented Lagrangian algorithm, can be used to find initial state and time-dependent optimal control and trajectory functions. To demonstrate our approach, we provide an alternative computation of the optimal control laws derived in \cite{GF22} for OCPs arising from odd circulant games modelled by the replicator dynamics. We also apply our approximation scheme to a quadratic OCP in the context of Newtonian mechanics in order to demonstrate the broader applicability of our method. 
\end{enumerate}

The remaining of this paper is organized as follows: in \ref{sec: Related Literature}, an overview of the current state of the literature addressing the numerical solution of OCPs is presented. Additionally, a discussion on the role of machine learning methods in the context of numerical optimal control is developed. Finally, odd circulant games are briefly introduced -- a class of evolutionary games giving rise to the game of rock-paper-scissors, contextualizing the main OCP of interest in this paper. \ref{sec: Problem Construction} introduces the main OCP for which the approximation scheme of this paper is applied, as well as the process of translating the OCP in parametric form into an unconstrained optimization problem, suitable for use with automatic differentiation. \ref{sec: Computational Experiments} presents a series of results from computational experiments performed for the main OCP of interest and for an additional OCP in the context of Newtonian mechanics. \ref{sec: Error Bounds} presents the derivation of MSE bounds for the multidimensional, truncated Fourier approximation considered in this work for the case of one and two-dimensional state spaces, inducing a general bound for problems with state space of dimension $n$. In light of the results obtained, \ref{sec: Conclusions} concludes this paper with an assessment of the approximation scheme suggested, along with future directions in overcoming some of the challenges and drawbacks of the proposed method.

\section{Related Literature and Preliminaries}\label{sec: Related Literature}
In this section, we provide a brief overview of the field of numerical optimal control. We also present an outline and discuss some recent results concerning the solution of OCPs under the modern lens of machine learning and function approximation. We then introduce the class of odd circulant evolutionary games, which we use as an example application.

\subsection{Numerical Solutions of Optimal Control Problems}
The theory of optimal control is one of the fundamental areas of modern engineering and physical sciences. Since its formal establishment in the 1950s by Richard Bellman, \cite{B52} followed by Lev Pontryagin \cite{P61}, the field has grown, taking on different names and serving even more diverse applications. The work in \cite{P21book} is a recent effort to create a better taxonomy and unified approach to this vast field. In practice, optimal control problems are solved numerically. We refer the reader to \cite{JTB20} for a consolidated text in the field.  In broader terms, these numerical methods are divided into two categories: direct and indirect methods. Direct methods translate the continuous problem into a discrete one, forming a nonlinear programming problem, which is then solved by standard open-source or commercial solver software. In their turn, indirect methods explore the necessary conditions from Pontryagin's Minimum Principle \cite{P61}, deriving a problem-related Hamiltonian and thus solving the system of differential equations arising therein, also known as the Euler-Lagrange equations \cite{R14,K17}. Both strategies rely on traditional numerical methods, such as those applied in the computation of solutions of constrained and unconstrained optimization problems, shooting methods for the solution of differential equations, and time-marching numerical integration techniques and algorithms for solving nonlinear systems of equations \cite{R14}. 
With the recent rise in interest in machine learning in all areas of applied sciences, a question is naturally raised when it comes to the development of the numerical OCP apparatus: what is the role of machine learning in leveraging these methods? Next, we present some recent developments that address this question.

\subsection{A Machine Learning Perspective on Optimal Control}
Work in the intersection of numerical methods for the solution of OCPs and machine learning has been of recent interest, involving both direct and indirect methods. The work in \cite{EP13, N17} proposes feedforward neural networks as surrogates for state, co-state and control in the computation of necessary conditions from the Hamiltonian system. In \cite{S18}, deep learning is also used, where feedforward neural networks are trained to learn the optimal control actions of an aircraft during landing. In the context of epidemiological control, the work in \cite{Y23} also approximates control functions by feedforward neural networks, using adjoint sensitivity analysis to solve an unconstrained optimization problem where the decision variables are composed of the parameters of the approximating neural network. Outside the scope of deep learning, \cite{N22, NG23} proposed an approximation of optimal controls as functions of both time and initial condition using Fourier series. 

In this paper, we adopt a machine learning perspective as well. As in a direct method, and by parameterizing both state and control by Fourier series of appropriate dimensions, we construct an unconstrained optimization problem to be solved with a gradient-requiring method, such as Conjugate Gradient Descent or Limited Memory Broyden–Fletcher–Goldfarb–Shanno (LBFGS). To overcome the complexity involved in the computation of the gradient required in the descent procedure, we exploit the availability of modern automatic differentiation software, also skewing from error-prone and difficult to implement numerical methods like finite differences \cite{B18,W64}. Furthermore, the use of a Fourier basis in the approximations places our approach under the umbrella of direct orthogonal collocation methods, sometimes also referred to as pseudospectral collocation. Typically, these methods involve orthogonal polynomials as so-called trial functions, like Lagrange \cite{EKR95, BHTR06} or Chebyshev \cite{VV88} polynomials, but work utilizing Fourier basis functions has also appeared in \cite{NY90, N22, Z22}. More specifically, we opt for a direct global orthogonal collocation \cite{HR08} approach. Different from a local approach, global collocation methods involve only one trial function on the time domain over which the function approximation is computed for a set of collocation points. Thus, for example, if we choose a Fourier basis to approximate a control function, this function will be defined over the entire time horizon $[0,T]$ of the OCP of interest. On the other hand, local methods partition the time domain into different segments each involving a trial function itself, thus requiring a set of parameters (Fourier or polynomial coefficients, for instance) for each segment. In this case, each segment has its own set of collocation points as well \cite{R14}. 

Additionally, in this work, we use multidimensional Fourier basis functions for the approximation of control and state functions over $[0,T]$, stepping into the realm of global orthogonal collocation methods. This choice also makes our approach appealing from the machine learning perspective, by combining optimization and automatic differentiation to find the best set of parameters making up the approximating functions. Another important aspect of the scheme we propose is that our control and state approximations are functions of both time and initial state condition, making it again suitable to a learning context and contrasting with the majority of work in numerical optimal control, which is primarily concerned with a single trajectory departing from a fixed initial condition. By solving the derived optimization problem given a finite set of initial conditions, we construct approximating optimal control and state surfaces (or hypersurfaces), that are approximately optimal for any initial condition in the range defined by the training set. This paper intends to generalize the work in \cite{N22, NG23}, which, to the authors' best knowledge, is the first of their kind in attempting to construct initial condition and time-dependent, continuous optimal control approximations.

A practical justification for the choice of basis function in our work is our intended application. As we present next, an OCP arises from a controlled rock-paper-scissors game dynamics, for which the optimal control is known to present oscillatory behaviour \cite{GF22}. Thus, we explore this aspect of our main test bed problem, justifying the use of Fourier approximations. However, as it will be demonstrated, this approximation scheme is also suitable to other types of OCPs. Next, we introduce the reader to the class of odd-circulant evolutionary games that shall give rise to one of the OCPs for which this paper proposes an approximated solution.

\subsection{Odd Circulant Evolutionary Games}
In the literature on classical and evolutionary game theory \cite{M12, W97}, odd-circulant games are defined as games whose payoff matrix posses a circulant form. These games are important in characterizing phenomena involving cyclic competition, such as in the famous children's game rock-paper-scissors \cite{G16}, where rock beats scissors, scissors beats paper and paper beats rock, in a cyclic form. Circulant matrices \cite{D79} are fully defined by their first row, and have the form, 
\begin{equation*}
\mathbf{A} =
\begin{pmatrix}
a_0 & a_{n-1} & a_{n-2} & \dots & a_1\\
a_1 & a_0 & a_{n-1} & \dots & a_2\\
\vdots & \vdots & \vdots & \ddots\\
a_{n-1} & a_{n-2} & a_{n-3} & \dots & a_0
\end{pmatrix}.
\end{equation*}
In $\mathbf{A}$, the second row is a cyclic permutation of the first row, and so it is with all the other rows. Take once again the example of rock-paper-scissors. Its payoff matrix is defined as,
\begin{equation*}
    \mathbf{L}_3 = \begin{pmatrix}
        0 & -1 &  1\\
        1 &  0 & -1\\
       -1 &  1 &  0\\
    \end{pmatrix},
\end{equation*}
where the cyclic structure can be easily observed. Now, we define the (also circulant) matrix $\mathbf{M}_3$, which serves as an actuating matrix determining the stability of the interior fixed point of this cyclic game,
\begin{equation*}
    \mathbf{M}_3 = 
    \begin{pmatrix}
        0 & 0 & 1\\
        1 & 0 & 0\\
        0 & 1 & 0
    \end{pmatrix}.
\end{equation*}
In a rock-paper-scissors game which is to be controlled toward equilibrium by an external control input $\gamma$ altering the payoff matrix structure, and thus targeting $\mathbf{u}^* = [\frac{1}{3}\;,\frac{1}{3}\;,\frac{1}{3}]$ in the unit simplex $\Delta_2 = \{\mathbf{u}\in\mathbb{R}^3\;:\;\mathbf{1}^\top\mathbf{u}=1,\;\mathbf{u}>\mathbf{0}\}$, a new payoff matrix arises,
\begin{equation}
    \mathbf{A}_3(\gamma) = \mathbf{L}_3 + \gamma\mathbf{M}_3.
    \label{eqn: payoff matrix}
\end{equation}
Here $\mathbf{u}$ denotes the vector containing the proportion of the total population playing each strategy. Notice how the matrix $\mathbf{M}_3$ has $1$'s where $\mathbf{L}_3$ also has $1$'s and zeros everywhere else. By making $\gamma \ge -1$ the game control input, we can choose an input $\gamma$ such that the payoff matrix $\mathbf{A}_3$ will naturally lead the system toward $\mathbf{u}^*$ \cite{GF22}. These definitions can be naturally generalized for circulant games with $N=5,7,9\dots$ strategies. An example is the rock-paper-scissors-Spock-lizard game \footnote{Developed by Sam Kass, see \url{http://www.samkass.com/theories/RPSSL.html}}. The dynamics of $\mathbf{u}$ under the action of $\gamma$ is described by the differential equation,
\begin{equation}
    \begin{aligned}
    &\dot{\mathbf{u}}_i = \mathbf{u}_i\left(\left(\mathbf{e}_i - \mathbf{u}\right)^\top\mathbf{A}_N(\gamma)\mathbf{u}\right)\\
    &\mathbf{u}(0) = \mathbf{u}_0.
    \end{aligned}
\end{equation}
As long as the initial condition $\mathbf{u}_0$ lies within the unit simplex $\Delta_{N-1}$, $\mathbf{u}(t)$ must also stay in $\Delta_{N-1}$ throughout the trajectory, regardless of the presence of $\gamma$ \cite{HS98}. Using \cref{eqn: payoff matrix},  we can write,
\begin{equation}
    \begin{aligned}
        &\dot{\mathbf{u}} = \mathbf{F}(\mathbf{u}) + \gamma\mathbf{G}(u)\\
        &\mathbf{u}(0) = \mathbf{u}_0,
        \label{eqn: Replicator Dynamics}
    \end{aligned}
\end{equation}
where we define $\mathbf{F}$, $\mathbf{G}: \Delta_{N-1} \rightarrow \mathbb{R}^n$ componentwise as,
\begin{equation}
    \begin{aligned}
        &\mathbf{F}_i(\mathbf{u}) = u_i\left(\left(\mathbf{e}_i - \mathbf{u}\right)^\top\mathbf{L}_N\mathbf{u}\right)\\
        &\mathbf{G}_i(\mathbf{u}) = u_i\left(\left(\mathbf{e}_i - \mathbf{u}\right)^\top\mathbf{M}_N\mathbf{u}\right).
        \label{eqn: F and G}
    \end{aligned}
\end{equation}
\cref{eqn: Replicator Dynamics} is the dynamics constraining the OCP introduced in the next section. 

\section{Problem Construction}\label{sec: Problem Construction}
Consider the OCP described in \cite{GF22},
\begin{equation}
    \begin{aligned}
    \min &\int_0^T \frac{1}{2}\lVert \mathbf{u} - \mathbf{u}^* \rVert^2 + \frac{r}{2}\gamma^2 \; dt \\
    s.t. \quad &\dot{\mathbf{u}} = \mathbf{F(u)} + \gamma\mathbf{G(u)} \\
    &\mathbf{u}(0) = \mathbf{u_0}, 
    \end{aligned}
    \label{eqn: mainOCP}
\end{equation}
where the state variable $\mathbf{u}, \mathbf{u}^* \in \mathbb{R}^d$ with control $\gamma$. Furthermore, $\mathbf{F}, \mathbf{G}$ are defined as in \cref{eqn: F and G} and $r$ is a scalar. Next, we propose an approximation for the control function $\gamma$.

\subsection{Approximation for Fixed Initial Condition}
Consider the OCP~(\ref{eqn: mainOCP}) with given $\mathbf{u}(0) = \mathbf{u_0}$. We introduce a controller approximation of the form,
\begin{equation}
    \gamma \approx \hat{\gamma}(t, \mathbf{a}, \mathbf{b}) = \sum_{m=1}^M a_m \sin\left(\frac{m \pi t}{T}\right) + \sum_{n=0}^N b_n \cos\left(\frac{n \pi t}{T}\right).
    \label{eqn: controlapprox}
\end{equation}
Similarly, a state approximation is made by,
\begin{equation}
    \mathbf{u} \approx \hat{\mathbf{u}}(t, \mathbf{C}, \mathbf{D}) = \sum_{m=1}^M \mathbf{c_m} \sin\left(\frac{m \pi t}{T}\right) + \sum_{n=0}^N \mathbf{d_n} \cos\left(\frac{n \pi t}{T}\right).
    \label{eqn: state approx}
\end{equation}
Now, the OCP~(\ref{eqn: mainOCP}) can be reformulated as,
\begin{equation}
    \begin{aligned}
    \min &\int_0^T \frac{1}{2}\lVert \mathbf{u} - \mathbf{u}^* \rVert^2 + \frac{r}{2}\gamma^2 \; dt \\
    s.t. \quad &\dot{\mathbf{u}} = \mathbf{F(u)} + \gamma\mathbf{G(u)} \\
    &\gamma(t) \approx \sum_{m=1}^M a_m \sin\left(\frac{m \pi t}{T}\right) + \sum_{n=0}^N b_n \cos\left(\frac{n \pi t}{T}\right)\\
    &\mathbf{u}(t) \approx \sum_{m=1}^M \mathbf{c_m} \sin\left(\frac{m \pi t}{T}\right) + \sum_{n=0}^N \mathbf{d_n} \cos\left(\frac{n \pi t}{T}\right)\\
    &\mathbf{u}(0) = \mathbf{u_0}.
    \end{aligned}
    \label{eqn: parOCP}
\end{equation}
The OCP~(\ref{eqn: parOCP}) is the parameterized version of OCP~(\ref{eqn: mainOCP}), and we shall denote it by $J(\mathbf{a},\mathbf{b}, \mathbf{C}, \mathbf{D}, \mathbf{u_0})$, where $\mathbf{a}$ and $\mathbf{b}$ are the vectors of Fourier coefficients in the approximation of the control variable $\gamma(t)$, with $\mathbf{a} = \left[a_1, \dots, a_M \right]^\top$ and $\mathbf{b} = \left[b_0, b_1, \dots, b_N \right]^\top$. In the state approximation, $\mathbf{C} = \left[\mathbf{c_1}, \mathbf{c_2}, \dots, \mathbf{c_M} \right]$ where $\mathbf{c_m} = \left[c_m^1, \dots, c_m^d \right]^\top$, thus $\mathbf{C}$ is a matrix with each column representing the set consisting of the $m^{th}$ Fourier coefficients in the approximation of the state space of dimension $d$, and $m=1,\dots,M$. Similarly, $\mathbf{D} = \left[\mathbf{d_1}, \mathbf{d_2}, \dots, \mathbf{d_M} \right]$ where $\mathbf{d_m} = \left[d_m^1, \dots, d_m^d \right]^\top$. Thus, the time continuous problem, \cref{eqn: mainOCP}, is cast into a nonlinear programming problem with a finite number of decision variables,
\begin{equation}
    \min_{\mathbf{a, b, C, D}} J\left(\mathbf{a, b, C, D, u_0}\right).
    \label{eqn: Discrete Optimization Problem - fixed x0}
\end{equation}
In this notation we make implicit the minimization of the parameterized objective functional $\int_0^T \frac{1}{2}\lVert \mathbf{u} - \mathbf{u}^* \rVert^2 + \frac{r}{2}\gamma^2 dt $ such that all the constraints in OCP~(\ref{eqn: parOCP}) are obeyed. Given the non-separability of the state dynamics $\dot{\mathbf{u}} = \mathbf{F(u)} + \gamma\mathbf{G(u)}$, we cannot express $\mathbf{u}$ as an approximation parameterized by $\mathbf{a}$ and $\mathbf{b}$ that would naturally arise from the integration of $\dot{\mathbf{u}}$ in OCP~(\ref{eqn: parOCP}), as done in \cite{N22, NG23}, in which only a Fourier approximation for the control is considered. Therefore, our suggested approximation of both state and control variables in this present work is justified. 
To simplify the expressions containing the proposed Fourier approximations, we recall that Fourier series can be written in the complex exponential form by using the following relations derived from Euler's formula, $\sin{\theta} + i\cos{\theta} = e^{(i\theta)}$,
\begin{equation}
    \cos{\theta} = \frac{e^{(i\theta)} + e^{(-i\theta)}}{2}, \quad \sin{\theta} = \frac{e^{(i\theta)} - e^{(-i\theta)}}{2i},
    \label{eqn: Euler relations}
\end{equation}
yielding the alternative expression for \cref{eqn: controlapprox} \cite{T76},
\begin{equation}
    \gamma \approx \hat{\gamma}(t, \mathbf{c}) = \sum_{k=-K}^K c_k e^{\frac{\pi i k t}{T}},
    \label{eqn: complex-exp Fourier}
\end{equation}
where $\mathbf{c} = \left[c_1, \dots, c_k\right]^\top$ and $i$ denotes the imaginary unit. Likewise, \cref{eqn: state approx} can be written in the complex exponential form. When the number of terms $2K+1$ is of particular importance, we shall denote $\hat{\gamma}(t, \mathbf{c})$ by $\hat{\gamma}_K(t)$. Additionally, to shorten the expressions involving trigonometric terms, we adopt the following notation through the remaining of this paper: for an argument of the type $\frac{m \pi t}{T}$ or $\frac{n \pi u_0}{\mathcal{U}_0}$ with $m, n \in \mathbb{Z_{+}}$, we make $\theta_m = \frac{m \pi t}{T}$ or $\phi_n = \frac{n \pi u_0}{\mathcal{U}_0}$, respectively. This notation use will become handy when representing approximations involving many basis functions. Next, we develop $2$-dimensional and $(d+1)-$dimensional Fourier approximations. 

\subsection{Approximation for a Range of Initial Conditions}
A solution to Problem~(\ref{eqn: Discrete Optimization Problem - fixed x0}) provides a set of coefficients such that \cref{eqn: controlapprox,eqn: state approx} are approximate optimal solutions to the OCP~(\ref{eqn: mainOCP}) with fixed initial condition $\mathbf{u_0}$. Now, consider the one dimensional state case $\mathbf{u} \in \mathbb{R}$, so $\mathbf{u_0} \in \mathbb{R}$. Given a range of feasible initial states of OCP~(\ref{eqn: mainOCP}) belonging to a given ordered set $U_0 =\{ u_0^1, u_0^2, \dots, u_0^j \} \subset \mathbb{R}$, and extending the Fourier approximation to the multivariable case, we propose the following approximations for the control and state trajectories of OCP~(\ref{eqn: mainOCP}),
\begin{multline}
\gamma \approx \hat{\gamma}(t, u_0, \mathbf{A, B, C, D}) = \\
\sum_{m,n=0}^{M,N } \bigl[ a_{mn}\sin\left(\theta_m\right)\sin\left(\phi_n\right) + b_{mn}\sin\left(\theta_m\right)\cos\left(\phi_n\right)\\
+c_{mn}\cos\left(\theta_m\right)\sin\left(\phi_n\right) + d_{mn}\cos\left(\theta_m\right)\cos\left(\phi_n\right) \bigl], 
    \label{eqn: control approx 2D}
\end{multline}
\begin{multline}
    \mathbf{u}_i \approx \hat{\mathbf{u}}_i(t, u_0, \mathbf{A}_i, \mathbf{B}_i, \mathbf{C}_i, \mathbf{D}_i) = \\
    \sum_{m,n=0}^{M,N} \bigl[ a^{u_i}_{mn}\sin\left(\theta_m\right)\sin\left(\phi_n\right) + b^{u_i}_{mn}\sin\left(\theta_m\right)\cos\left(\phi_n\right)\\
+c^{u_i}_{mn}\cos\left(\theta_m\right)\sin\left(\phi_n\right) + d^{u_i}_{mn}\cos\left(\theta_m\right)\cos\left(\phi_n\right) \bigl],
    \label{eqn: state approx 2D}
\end{multline}
for $i=1,\dots,d$, where we organize the coefficients of $\hat{\gamma}$ as matrices $\mathbf{A, B, C}$ and $\mathbf{D}$, of dimension $\mathbb{R}^{(M+1)\times (N+1)}$, and the coefficients of $\hat{u}$ as $(M+1)\times(N+1)\times d-$dimensional arrays. To make the notation more compact, we shall refer to $\mathbf{\Theta}$ and $\mathbf{\Theta}_u$ as the collection of coefficients parameterizing $\hat{\gamma}$ and $\hat{u}$, respectively. Additionally, $\mathcal{U}_0 = u_0^j - u_0^1$.  
In such a case, where the aforementioned approximations are functions of both time and initial condition, the parameterized OCP of interest becomes,
\begin{equation}
    \begin{aligned}
    \min &\int_0^T \frac{1}{2}\lVert \mathbf{u} - \mathbf{u}^* \rVert^2 + \frac{r}{2}\gamma^2 \; dt \\
    s.t. \quad &\dot{\mathbf{u}} = \mathbf{F(u)} + \gamma\mathbf{G(u)} \\
    &\gamma(t, u_0, \mathbf{\Theta}) \approx \hat{\gamma}(t, u_0, \mathbf{\Theta})\\
    &\mathbf{u}_i(t, u_0, \mathbf{\Theta}_u) \approx \hat{\mathbf{u}}_i(t, u_0, \mathbf{\Theta}_u)\\
    &\mathbf{u}(0) \in U_0 \subset \mathbb{R}^2,
    \end{aligned}
    \label{eqn: parOCP 2D}
\end{equation}
which we aim to approximately solve by formulating an optimization problem equivalent to that of OCP~(\ref{eqn: Discrete Optimization Problem - fixed x0}), yet now considering the set of initial conditions in $U_0$,
\begin{equation}
    \min_{\mathbf{\Theta}, \mathbf{\Theta}_u} \sum_{u_0 \in U_0}  J\left(\mathbf{\Theta}, \mathbf{\Theta}_u\right)
    \label{eqn: Discrete Optimization Problem - various x0}.
\end{equation}
A solution to Problem~(\ref{eqn: Discrete Optimization Problem - various x0}) allows us to construct an approximate optimal control surface,
\begin{equation}
\{ \bm{\gamma}(t, u_0, \bm{\Theta}, \bm{\Theta}_u) \;\vert\; (t,u_0) \in \left[0,T\right] \times \left[{u_0}^1, {u_0}^j \right] \}
\end{equation}
for OCP~(\ref{eqn: parOCP 2D}). 

Finally, consider the general case $\mathbf{u} \in \mathbb{R}^d$. If we want a $(d+1)-$dimensional approximation of the control $\gamma(t, \mathbf{u}_0)$ where $\mathbf{u}_0 \in \mathbb{R}^d$, we make,
\begin{multline}
        \gamma(t, \mathbf{u_0}, \mathbf{\Theta}, \mathbf{\Theta}_u) \approx\\
        \sum^{M, N_1, N_2, \dots, N_d}_{m, n_1, n_2, \dots, n_d = 0}  \bigl[ a^1_{m, n_1, n_2,\dots,n_d}\sin{(\theta_m)}\sin{(\phi_{m_1})} \dots \sin{(\phi_{m_d})} + \dots +\\  a^{j}_{m, n_1, n_2,\dots,n_d}\cos{(\theta_m)}\cos{(\phi_{m_1})} \dots     
  \cos{(\phi_{m_d})}\bigr],
\end{multline}
where $j=2^{(d+1)}$ (when all the $2^{(d+1)}$ basis functions are used) and, for computational purposes, $\mathbf{a}^i$ is the $(d+1)-$dimensional array containing all the coefficients $a^i_{m, n_1,\dots,n_d}$, for $i=1,\dots,2^{(d+1)}$. 

\subsection{Optimization Scheme}
The nonseparability of the dynamics of OCP~(\ref{eqn: mainOCP}) led us to introduce two sets of parameters in the approximation of control and state by $\hat{\gamma}$ and $\mathbf{\hat{u}}$, respectively. Therefore, the priced-out dynamics must be present in the parameterized objective functional. In the literature of nonlinear programming, there are many known techniques to tackle this translated unconstrained (parameterized) optimization problem. In this work, we opt to use the Augmented Lagrangian approach as presented in \cite{bazaraa2013nonlinear}. Among many other reasons, the Augmented Lagrangian skews from computational difficulties associated with ill-conditioning and ever-increasing multiplier parameters, doing so by leveraging the positive aspects of both penalty functions and Lagrange multiplier approaches. Nonetheless, the optimization approach described in \cite{bazaraa2013nonlinear} is suited for a static nonlinear optimization problem, forcing us to account for the inherent dynamics of OCP~(\ref{eqn: parOCP}). 

First, consider the (static) optimization problem,
\begin{equation}
    \min_{x} \; f(\mathbf{x}) \quad \text{s.t.} \quad \mathbf{x} \in D = \{\mathbf{x} \in \mathbb{R}^n\;\vert\;\mathbf{h}(\mathbf{x}) = 0\},
    \label{eqn: ECNLP}
\end{equation}
where $f(\mathbf{x}): \mathbb{R}^n \longrightarrow \mathbb{R}$ and $\mathbf{h}(\mathbf{x}): \mathbb{R}^n \longrightarrow \mathbb{R}^l$. The Augmented Lagrangian function is defined as,
\begin{equation}
    \mathcal{L}(\mathbf{x}, \boldsymbol{\upsilon}; \boldsymbol{\mu}) = f(\mathbf{x}) + \sum_{i=1}^l \upsilon_i h_i(\mathbf{x}) + \sum_{i=1}^l \mu_i h_i^2(\mathbf{x}),
    \label{eqn: Lagrangian}
\end{equation}
where $\mathbf{x}$ and the Lagrange multipliers $\bm{\upsilon} = [\upsilon_1 \dots \upsilon_l]^\top$ are decision variables and $\boldsymbol{\mu} = [\mu_1 \dots \mu_l]^\top$ are penalty parameters. It can be shown \cite{izmailov2007otimizaccao} that a stationary point of the constrained problem \cref{eqn: ECNLP} is also a stationary point for the Augmented Lagrangian \cref{eqn: Lagrangian}. 

Now, in light of the OCP~(\ref{eqn: parOCP}), we define the functions,
\begin{align}
    &f(\mathbf{\Theta}, \mathbf{\Theta}_u) = \int_0^T \frac{1}{2}\lVert \hat{\mathbf{u}}(\mathbf{\Theta}_u) - \mathbf{u}^* \rVert^2 + \frac{r}{2}\hat{\gamma}^2(\mathbf{\Theta}) \; dt \label{eqn:fpar}\\
    &h_{1}(\mathbf{\Theta}, \mathbf{\Theta}_u) = \int_0^T \frac{1}{2}\left\lVert\frac{d}{dt}{\hat{u}}_i(\mathbf{\Theta}_u) - F_i(\hat{u}(\mathbf{\Theta}_u)) - \hat{\gamma}(\mathbf{\Theta})G_i(\hat{u}(\mathbf{\Theta}_u))\right\rVert^2\;dt\label{eqn:h1par}\\
    &h_{i,2}(\mathbf{\Theta}, \mathbf{\Theta}_u) = \int_0^T \max\left(0, -\hat{u}_i(\mathbf{\Theta}_u)\right)\;dt\label{eqn:h2par}\\
    &h_{3}(\mathbf{\Theta}, \mathbf{\Theta}_u) = \int_0^T \frac{1}{2}\lVert \mathbf{1}^\top \hat{\mathbf{u}}(\mathbf{\Theta}_u) - 1 \rVert^2\;dt\label{eqn:h3par}\\
    &h_{4}(\mathbf{\Theta}, \mathbf{\Theta}_u) = \frac{1}{2}\lVert \mathbf{u}_0 - \hat{\mathbf{u}}(0, \mathbf{\Theta}_u)\rVert^2 \; dt,\label{eqn:h4par}
\end{align}
where \cref{eqn:h2par} is defined for all $i \in \{1,\dots,d\}$. \cref{eqn:fpar} is the parameterized objective functional of the original OCP~(\ref{eqn: mainOCP}), wherein $t$ is integrated out, thus making it a function of $\mathbf{\Theta}$ and $\mathbf{\Theta}_u$ only.  \cref{eqn:h1par}, by definition, is the distance in functional space, that is, the $\ell_2-$norm between the time derivative, 
\begin{equation*}
\frac{d}{dt}\hat{\mathbf{u}}(\mathbf{\Theta}_u)
\end{equation*}
and the also parameterized enforced dynamics, arising from the equality constraint in OCP~(\ref{eqn: mainOCP}). In words, we want to force the evolution of the parameterized state trajectory to match the right-hand side of the dynamics in OCP~(\ref{eqn: mainOCP}) $\forall\;t \in [0,T]$, accomplishing this by making $h_{1}(\mathbf{\Theta}, \mathbf{\Theta}_u)=0$. \cref{eqn:h3par} forces the state space to lie in the unit simplex $\Delta_{d-1}$ and \cref{eqn:h4par} forces the trajectory to start at $\mathbf{u}_0$. Substituting \cref{eqn:fpar,eqn:h1par,eqn:h2par,eqn:h3par,eqn:h4par} into the Lagrangian, \cref{eqn: Lagrangian}, (with decision variables $\mathbf{\Theta}$ and $\mathbf{\Theta}_u$), we obtain the following expression for $\mathcal{L}(\mathbf{\Theta}, \mathbf{\Theta}_u; \boldsymbol{\upsilon}, \boldsymbol{\mu})$,
\begin{multline}
f(\mathbf{\Theta}, \mathbf{\Theta}_u) + \upsilon_{1} h_{1}(\mathbf{\Theta}, \mathbf{\Theta}_u) +\\
\sum_{i=1}^d \left(\upsilon_{i,2} h_{i,2}(\mathbf{\Theta}, \mathbf{\Theta}_u)\right) +\upsilon_{3} h_{3}(\mathbf{\Theta}, \mathbf{\Theta}_u)  + \upsilon_{4} h_{4}(\mathbf{\Theta}, \mathbf{\Theta}_u)
        + \mu_{1} h_{1}^2(\mathbf{\Theta}, \mathbf{\Theta}_u) +\\
        \sum_{i=1}^d \left(\mu_{i,2} h_{i,2}^2(\mathbf{\Theta}, \mathbf{\Theta}_u) \right) +\mu_{3} h_{3}^2(\mathbf{\Theta}, \mathbf{\Theta}_u)  + \mu_{4} h_{4}^2(\mathbf{\Theta}, \mathbf{\Theta}_u),
    \label{eqn: Lagrangian2}
\end{multline}
where $\boldsymbol{\upsilon} = \left[{\upsilon}_1 \dots {\upsilon}_4\right]$ and $\boldsymbol{\upsilon} = [{\mu}_{1}\dots{\mu}_{4}]$.
Next, we illustrate how the optimization procedure is performed via the Augmented Lagrangian, shown in \cref{alg: AL}. For exposition purposes, we depict it with a naive Gradient Descent, but in practice we employ other first-order descent methods, like LBFGS and Conjugate Gradient. More on the specific method used in this work will be discussed in \cref{sec: Computational Experiments}. It is also important to note that even though \cref{alg: AL} is presented with \cref{eqn: Lagrangian2} in mind (and therefore the objective functional and constraints from OCP~(\ref{eqn: parOCP})), this algorithm can be applied to any other continuous-time OCP of the same nature, as long as all constraints are put into equality form, as done with the constraints in  \cref{eqn:h2par}. This feature will be demonstrated in \cref{sec: Computational Experiments}, where we apply \cref{alg: AL} for a quadratic control OCP. The general structure of \cref{alg: AL} is inspired by the Augmented Lagrangian algorithm presented in \cite{bazaraa2013nonlinear}.   

\begin{algorithm}
\caption{Minimization of $\mathcal{L}(\mathbf{\Theta}, \mathbf{\Theta}_u; \boldsymbol{\upsilon}, \boldsymbol{\mu})$: Augmented Lagrangian}\label{alg: AL}
\begin{algorithmic}
\Require Initialize the coefficients $\mathbf{\Theta}^0, \mathbf{\Theta}_u^0$, initial control $\hat{\gamma}^0 = \hat{\gamma}(\mathbf{\Theta}^0)$ and state $\hat{\mathbf{u}}^0 = \hat{\mathbf{u}}(\mathbf{\Theta}^0_u)$, gradient norm tolerance $\epsilon>0$, step length $\alpha>0$, outer loop iterations limit $\ell_{lim}$, violation tolerance $\tau>0$, violation reduction factor $0 < \lambda \le 1$, Lagrange multipliers $\boldsymbol{\upsilon}^0 > 0$, Lagrange multiplier update factor $\lambda_j > 0,\; j=1,\dots,4$, penalty terms $\boldsymbol{\mu}^0 > 0$ and penalty term update factor $\lambda^j > 0,\;j=1,\dots,4$.  
\State $\ell \gets 0$;
\State Set initial violation to $\nu^0 = \infty$;
\While{$\nu^\ell = \max\{ h_{1}(\mathbf{\Theta}^\ell, \mathbf{\Theta}_u^\ell), h_{2}(\mathbf{\Theta}^\ell, \mathbf{\Theta}_u^\ell), h_{3}(\mathbf{\Theta}^\ell, \mathbf{\Theta}_u^\ell), h_{4}(\mathbf{\Theta}^\ell, \mathbf{\Theta}_u^\ell)\} > \tau$ and $\ell > \ell_{lim}$} 
    \State $k \gets 0$;
    \State Compute $\nabla \mathcal{L}(\mathbf{\Theta}^0, \mathbf{\Theta}_u^0; \boldsymbol{\upsilon}^0, \boldsymbol{\mu}^0)$;
\While{$\lVert\nabla \mathcal{L}(\mathbf{\Theta}^k, \mathbf{\Theta}_u^k; \boldsymbol{\upsilon}^\ell, \boldsymbol{\mu}^\ell)\rVert \ge \epsilon$}
    \State Update coefficients $\mathbf{\Theta}^k \gets \mathbf{\Theta}^k - \alpha \nabla_{\mathbf{\Theta}} \mathcal{L}(\mathbf{\Theta}^k, \mathbf{\Theta}_u^k; \boldsymbol{\upsilon}^\ell, \boldsymbol{\mu}^\ell)$;
    \State Update coefficients $\mathbf{\Theta}_u^k \gets \mathbf{\Theta}_u^k - \alpha \nabla_{\mathbf{\Theta}_u} \mathcal{L}(\mathbf{\Theta}^k, \mathbf{\Theta}_u^k; \boldsymbol{\upsilon}^\ell, \boldsymbol{\mu}^\ell)$;
    \State Update control approximation $\hat{\gamma}^k = \hat{\gamma}(\mathbf{\Theta}^k)$ and trajectory $\hat{\mathbf{u}}^k = \hat{\mathbf{u}}(\mathbf{\Theta}_u^k)$;
    \State Compute $\nabla \mathcal{L}(\mathbf{\Theta}^k, \mathbf{\Theta}_u^k; \boldsymbol{\upsilon}^\ell, \boldsymbol{\mu}^\ell)$;
    \State $k \gets k + 1$;
\EndWhile
    \If{$\nu^\ell < \lambda \nu^{\ell - 1}$}
        \State Update Lagrange multipliers $\upsilon_j \gets \lambda_j \upsilon_j h_j(\mathbf{\Theta}^k, \mathbf{\Theta}_u^k), \; j=1,\dots,4$;
    \Else
        \State For each $h_j(\mathbf{\Theta}^k, \mathbf{\Theta}_u^k) \ge \lambda \nu^{\ell - 1}, \; j=1,\dots,4$, update penalty term $\mu_j \gets \lambda^j\mu_j$;
    \EndIf
\EndWhile
\end{algorithmic}
\end{algorithm}

\subsection{Automatic Differentiation}
In order to perform a first order descent procedure of Problem~(\ref{eqn: Lagrangian}) to find the approximated optimal control, as in \cite{N22}, the derivatives of the Augmented Lagrangian with respect to each Fourier coefficient in the approximation of control and state variables need to be computed, either approximately or exactly. In the former case, one could opt for a finite differences scheme or for an automatic differentiation (forward or backward) computation. In the latter, the analytical derivatives need to be derived. Given the intricacy of the expressions comprising both the original objective function $f(\cdot)$ and the constraints $h(\cdot)$, their closed-form, analytical derivatives become intractable, especially as the dimensionality of the OCP grows. Therefore, for the remainder of this paper, we compute the derivatives of $\mathcal{L}(\cdot)$ with respect to the Fourier coefficients through the use of automatic differentiation, a feature currently available for all major scientific and numerical computing software and programming languages. We note that derivative free methods can also be used, although this choice would skew our approach from a mainstream machine learning framework. 

\section{Computational Experiments and Results}\label{sec: Computational Experiments}
In this section, we implement the multivariable, Fourier-based control and state approximations discussed thus far for two distinct OCPs. The first problem involves a simple single particle whose dynamics obeys the laws of Newtonian mechanics. The optimal control for this problem can be found analytically, and its purpose is to illustrate the effectiveness of the proposed approximations. The second, and more complex case, is the OCP presented in \cite{GF22}, describing the evolutionary, odd-circulant game of rock-paper-scissors under the action of a centralizing agent (control) seeking to bring the strategies into equilibrium in the unit simplex. 

\subsection{Problem 1:  Linear Quadratic Control}
This problem is adapted from \cite{LVS12}. A particle of unit mass moves along a line with position $x(t)$ and velocity $\dot{x}(t)$, from time $t=0$ up to $t=T$, while being subject to an acceleration $\ddot{x}(t) = \gamma(t)$. The state space of this problem is defined by $\mathbf{x}(t) = \left[x(t) \; \dot{x}(t)\right]^{\top}$, with initial and final conditions $\mathbf{x}(0) = \left[x(0) \; \dot{x}(0)\right]^{\top}$and $\mathbf{x}(T) = \left[x(T) \; \dot{x}(T)\right]^{\top}$, respectively. Starting at $\mathbf{x}(0)$, we wish to bring this particle to $\mathbf{x}(T)$ by choosing an acceleration $\gamma(t) \in[0,T]$ in such a way that the value of $\int_0^T r\gamma(t)^2 \; dt$ is minimized, where $r$ is just a constant. In words, we want to bring the particle from initial to final condition with minimum effort. The optimal control formulation is,
\begin{equation}
    \begin{aligned}
        \min_\gamma &\int_0^T r\gamma(t)^2 \; dt\\
        &\dot{\mathbf{x}} = \begin{bmatrix} 0 & 1\\ 0 & 0 \end{bmatrix} \mathbf{x}(t) + \begin{bmatrix}
            0 \\ 1 \end{bmatrix} \gamma(t)\\
        &\mathbf{x}(0) = [x_0 \; v_0]^\top\\
        &\mathbf{x}(T) = [x_T \; v_T]^\top,
    \end{aligned}
    \label{eqn: toyproblem}
\end{equation}
with analytical solution,
\begin{equation}
\gamma^*(t) = \begin{bmatrix}
    \frac{6T - 12t}{T^3} & \frac{-2T + 6t}{T^2}
\end{bmatrix} \begin{bmatrix}
    \mathbf{x}(T) - \begin{bmatrix}
        1 & T\\ 0 & 1
    \end{bmatrix} \mathbf{x}(0)
\end{bmatrix}, 
\end{equation}
which shall be compared to the approximations to be computed. First, we consider a simpler case with fixed initial velocity and varying initial positions. 

\subsubsection{Varying Initial Position and Fixed Initial Velocity}
Fixing the initial velocity $\dot{x}(0) = 1$, we take from the set $X_0 = \{0, 0.5, \dots, 4.5, 5.0\}$ initial positions and approximate control and state as in \cref{eqn: control approx 2D} and \cref{eqn: state approx 2D}, respectively.  Then, we solve the unconstrained optimization problem with the Augmented Lagrangian as being defined as,
\begin{equation}
    \mathcal{L}(\Theta, \Theta_x) = \sum_{x_0 \in X_0}\left[ \int_0^T r\hat{\gamma}(t, x_0, \Theta)^2 \; dt \; \right] + \sum_{i=1}^3 \left[ \upsilon_i h_i + \mu_i h_i^2\right],
    \label{eqn: Lagrangian toy}
\end{equation}
where,
\begin{align}
        h_1(\Theta, \Theta_x) &= \sum_{x_0 \in X_0} \left[ \int_0^T \frac{1}{2} \left\lVert \frac{d}{dt}\hat{\mathbf{x}}(t, x_0, \Theta_x) - \mathbf{g}(t, x_0, \Theta, \Theta_x) \right\rVert^2 \; dt\right] \label{eqn: h1 toy}\\
        h_2(\Theta_x) &= \sum_{x_0 \in X_0} \left[\frac{1}{2} \lVert \hat{\mathbf{x}}(T, x_0, \Theta_x) - \mathbf{x}(T) \rVert^2 \right] \label{eqn: h2 toy}\\
        h_3(\Theta_x) &= \sum_{x_0 \in X_0} \left[\frac{1}{2} \lVert \hat{\mathbf{x}}(0, x_0, \Theta_x) - \mathbf{x}(0) \rVert^2 \right]. \label{eqn: h3 toy}
\end{align}
\cref{eqn: h1 toy} involves the time derivative of the state approximation, \textit{i.e.}, 
\begin{equation}
D_t\hat{\mathbf{x}}(t, x_0, \Theta_x) = \frac{d}{dt}[\hat{\mathbf{x}}(t, x_0, \Theta_x)],
\end{equation} 
which can be obtained directly from $\hat{\mathbf{x}}(t, x_0, \Theta_x)$, involving only trigonometric terms of easy computation. Furthermore, $\mathbf{g}(\cdot)$ is the (parameterized) right-hand side of the dynamics in OCP~(\ref{eqn: toyproblem}). Finally,  \cref{eqn: h1 toy} forces the time-derivative of $\mathbf{x}$ to equate to the right-hand side of the dynamics, and \cref{eqn: h2 toy} and \cref{eqn: h3 toy} force the initial and final conditions, respectively, over all initial positions in $X_0$. 
We solve this problem twice, each time picking a different number of Fourier coefficients $M + 1$ and $N+1$ in \cref{eqn: control approx 2D} and \cref{eqn: state approx 2D}. To simplify the implementation, we make $M=N$ throughout the remainder of this paper. 
In \cref{tab: Toy1}, we report, with respect to the approximate control, mean square error ($MSE_\gamma$), mean absolute percentage error ($MAPE_\gamma$) and  mean absolute error ($MAE_\gamma$). The number of iterations in the Augmented Lagrangian \cref{alg: AL} inner loop is, $k$ and that of the outer loop is denoted by $\ell$. The percentage error of the optimal cost $J^*$ over all initial conditions is denoted by, $J^*_{\%error}$ and we take $r=4$. Finally, for the experiments in \cref{tab: Toy1}, the choice of gradient requiring method was the Conjugate Gradient (CG) \cite{HZ06, HZ13} with the line search procedure implemented by the same authors and gradient norm tolerance $\epsilon$ and, the  LBFGS as in \cite{NW99}, using the same line search procedure as the CG alternative.  

\begin{table}[htbp]
    \centering
    \caption{Summary of computational results of Problem 1 with fixed $\dot{x}(0)$.}
    \label{tab: Toy1}
    \begin{tabular}{ccccccccc}
        \hline
        Method & $M=N$ & $k$ & $\ell$ & $MSE_\gamma$ & $MAPE_\gamma$ & $MAE_\gamma$ & $\epsilon$ & $J^{*}_{\%error}$ \\
        \hline
        \hline
        \multirow{4}{*}{CG} &$3$ & $1e+3$ & $20$ & $4.97e-2$ & $258.83\%$ & $1.93e-1$ & $10e-6$ & $23.76\%$ \\
        & $4$ & $1e+3$ & $30$ & $4.02e-4$ & $15.95\%$ &  $1.53e-2$  & $10e-6$  & $2.08\%$ \\
        & $4$ & $3e+3$ & $30$ & $4.91e-4$ & $18.12\%$ &  $1.70e-2$  & $10e-6$  & $2.07\%$ \\
        & $5$ & $1e+3$ & $15$ &  $1.56e-3$ & $17.31\%$ & $2.94e-2$ & $10e-4$ & $17.54\%$ \\ 
        \hline
        \multirow{2}{*}{LBFGS} & $4$ & $5e+3$ & $12$ & $4.40e-3$ & $33.96\%$ & $4.94e-2$ & $10e-4$ & $28.02\%$\\
        & $4$ & $3e+3$ & $30$ & $4.84e-4$ & $17.85\%$ & $1.70e-2$ & $10e-6$ & $1.78\%$\\
        \hline
    \end{tabular}
\end{table}

\begin{figure}[ht!]
\centering
    \centering
    \includegraphics[width=0.6\textwidth]{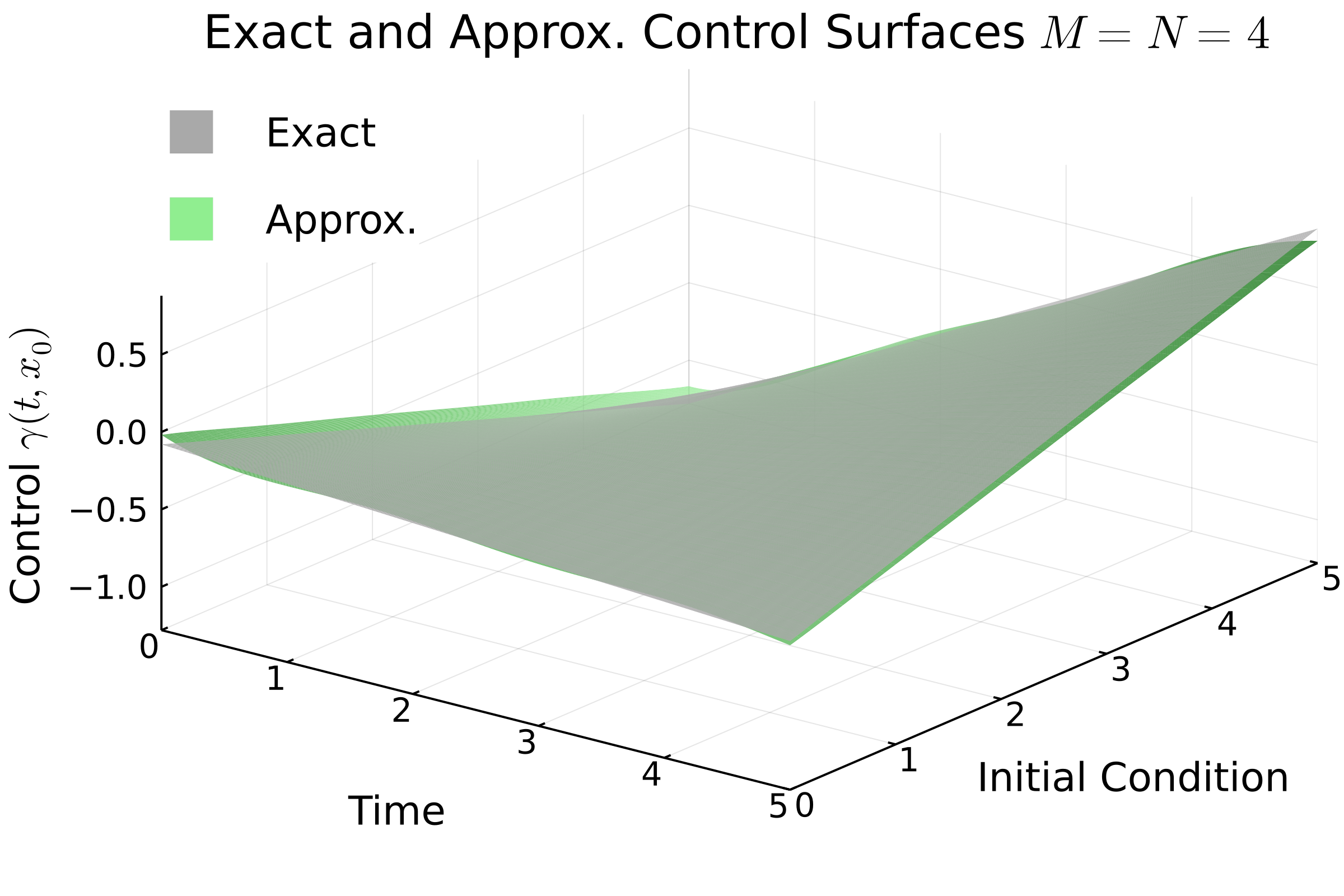}
     \qquad
    \includegraphics[width=0.6\textwidth]{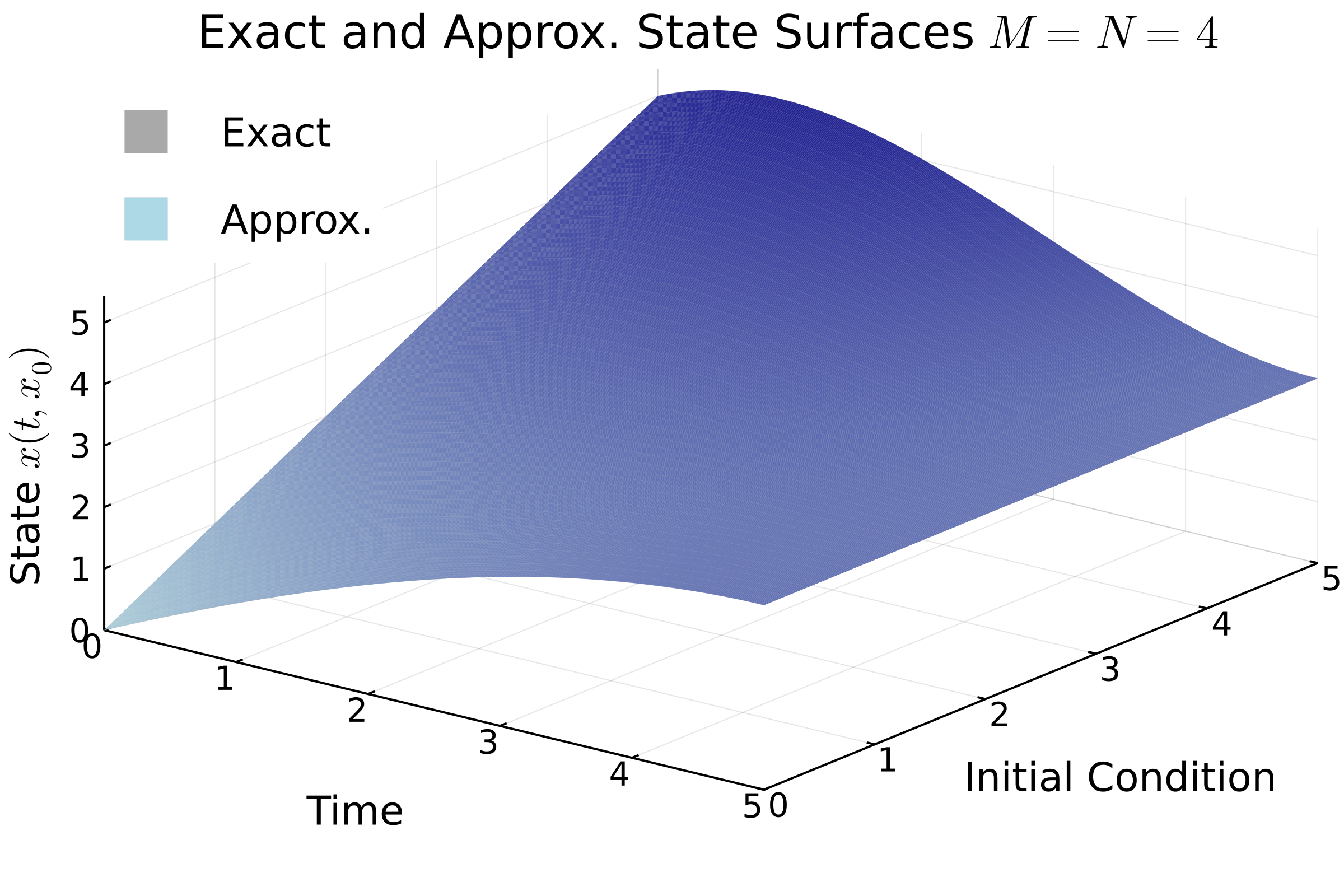}
    \caption{(Top) Control surface approximation for $M=N=4$. (Bottom) Corresponding optimal state (position) surface. Initial velocity fixed at $\dot{x}=1.0$.}
    \label{fig: toyMN4}
\end{figure}

\begin{figure}[ht!]
\centering
    \centering
    \includegraphics[width=0.6\textwidth]{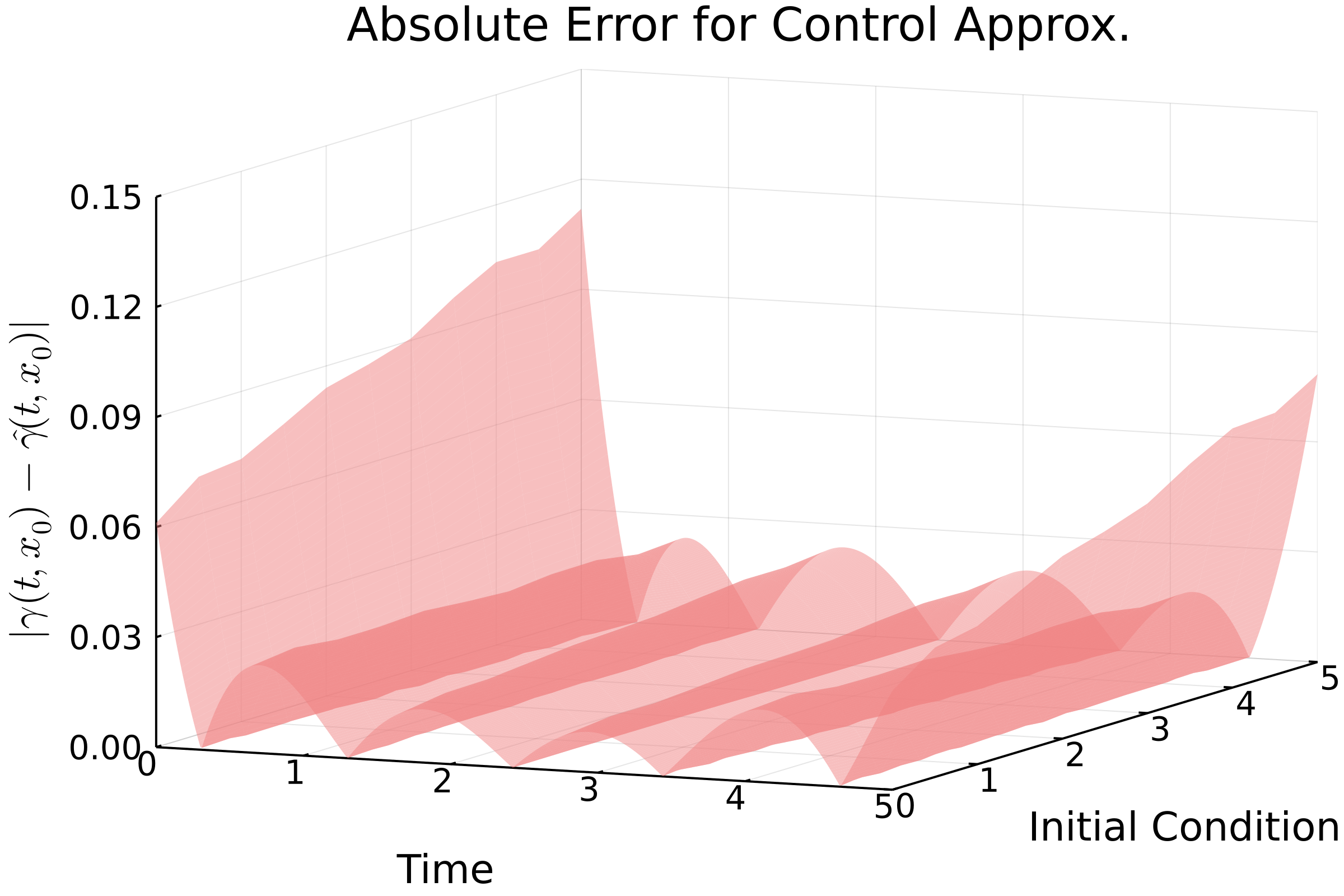}
    \qquad
    \includegraphics[width=0.6\textwidth]{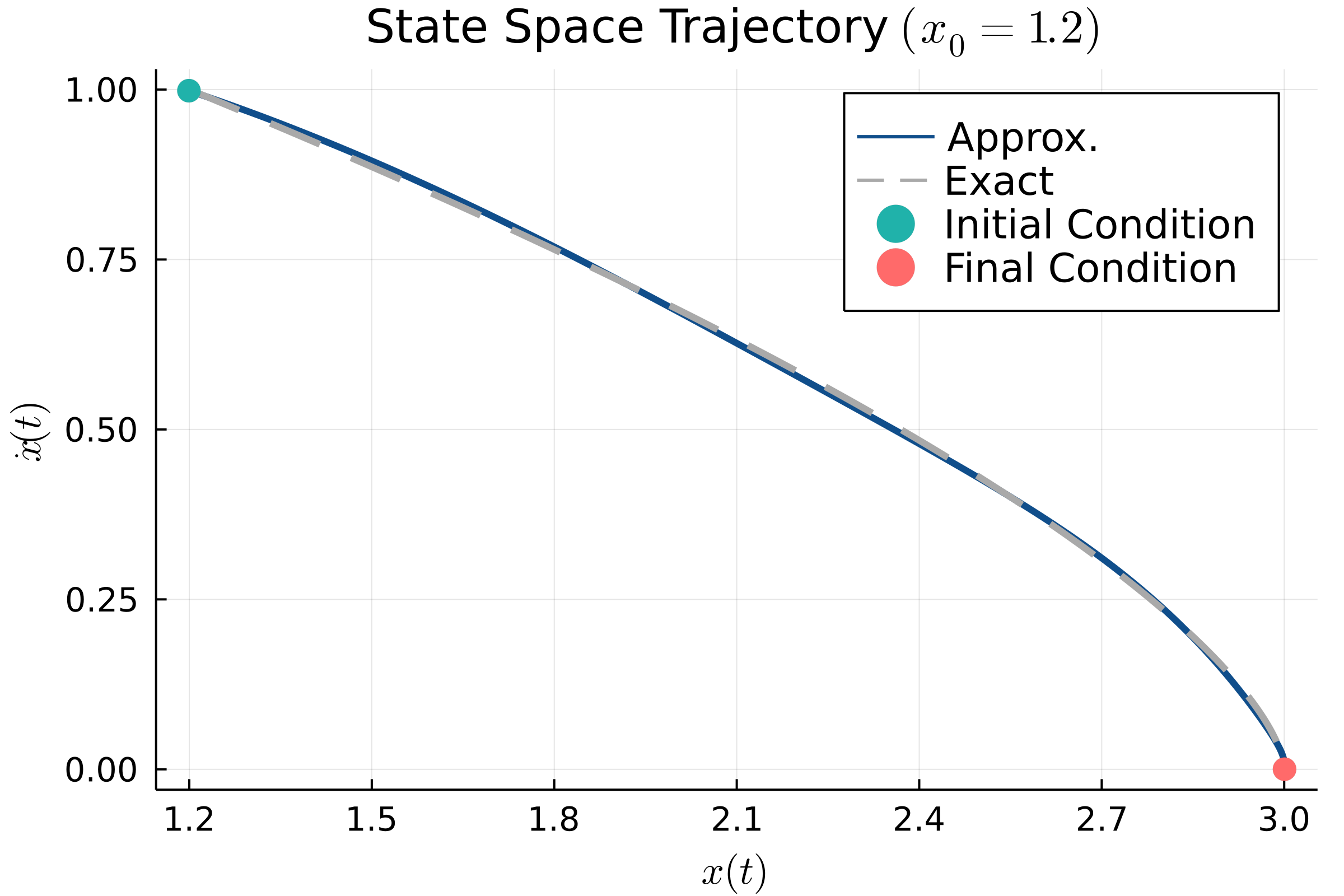}
    \caption{(Top) Absolute error for the control approximation in \cref{fig: toyMN4}. (Bottom) State space trajectory for the initial condition $x_0=1.2$ for the approximating control of \cref{fig: toyMN4}.}
    \label{fig: toyMN4 - 2}
\end{figure}

\cref{fig: toyMN4} shows the plots for both approximated optimal control and its corresponding optimal position trajectory surfaces. Additionally, \cref{fig: toyMN4 - 2} depicts the absolute error for the control surface of \cref{fig: toyMN4}, alongside a sample state space trajectory departing from an initial position $x_0 = 1.2 \notin X_0$. This demonstrates that our proposed control approximation cannot only be applied to trajectories starting at all $x_0 \in X_0$, but also to trajectories departing from initial positions that lie within the continuous set $[0, 5]$. 

\subsubsection{Varying Initial Position and Initial Velocity}
We now compute $3-$dimensional Fourier approximations of the form $\hat{\gamma}(t, x_0, \dot{x}_0, \Theta)$ and $\hat{\mathbf{x}}(t, x_0, \dot{x}_0, \Theta_x)$, where both variables of the state space have initial conditions $(x_0, \dot{x}_0) \in X_0 = \{0,1.0,\dots,5.0\}\times\{1.0,1.5,\dots,3.0\}$. The Augmented Lagrangian has the same form of \cref{eqn: Lagrangian toy}, but now summing over all $\mathbf{x}(0)$ in $X_0$. In this case, the optimal control and state surfaces lie in the $4-$dimensional Euclidean space, and cannot be visualized. To overcome this, we present the approximated optimal surfaces for different values of initial velocity, all derived from $\hat{\gamma}(t, x_0, \dot{x}_0, \Theta^*)$ and $\hat{\mathbf{x}}(t, x_0, \dot{x}_0, \Theta_x^*)$.  We report results for two different experiments, choosing 3 and 4 Fourier coefficients for each variable, that is, $M=N_1=N_2=2$ and  $M=N_1=N_2=3$. Such values are used because these experiments have shown to be computationally expensive, and this issue will be addressed in future work. Instances marked with an asterisk in \cref{tab: Toy2} indicate experiments where only half of the total number of basis functions were used in the truncated Fourier approximation. Both CG and LBFGS algorithms were employed. 
\begin{table}[htbp]
    \centering
    \caption{Summary of computational results of Problem 1 for varying $\dot{x}(0)$. The constant parameters used in all the experiments are $\ell=12$ and $\epsilon=10e-4$. }
    \label{tab: Toy2}
    \begin{tabular}{ccccccc}
        \hline
        Method & $M = N_i$ & $k$ & $MSE_\gamma$ & $MAPE_\gamma$ & $MAE_\gamma$ & $J^{*}_{\%error}$ \\
        \hline
        \hline
        \multirow{2}{*}{CG} & $2$ & $1e+3$ & $7.73e-2$ & $109.65\%$ & $1.96e-1$ & $22.73\%$ \\
        & $3^*$ & $1e+3$ & $1.91e-1$ & $271.54\%$ & $3.56e-1$ & $19.99\%$ \\
        \hline
        \multirow{1}{*}{LBFGS} & $2^*$ & $5e+3$ & $1.19e-1$ & $111.18\%$ & $2.45e+1$ & $29.05\%$\\
        
        \hline
    \end{tabular}
\end{table}

\begin{figure}[ht!]
\centering
    \centering
    \includegraphics[width=0.6\textwidth]{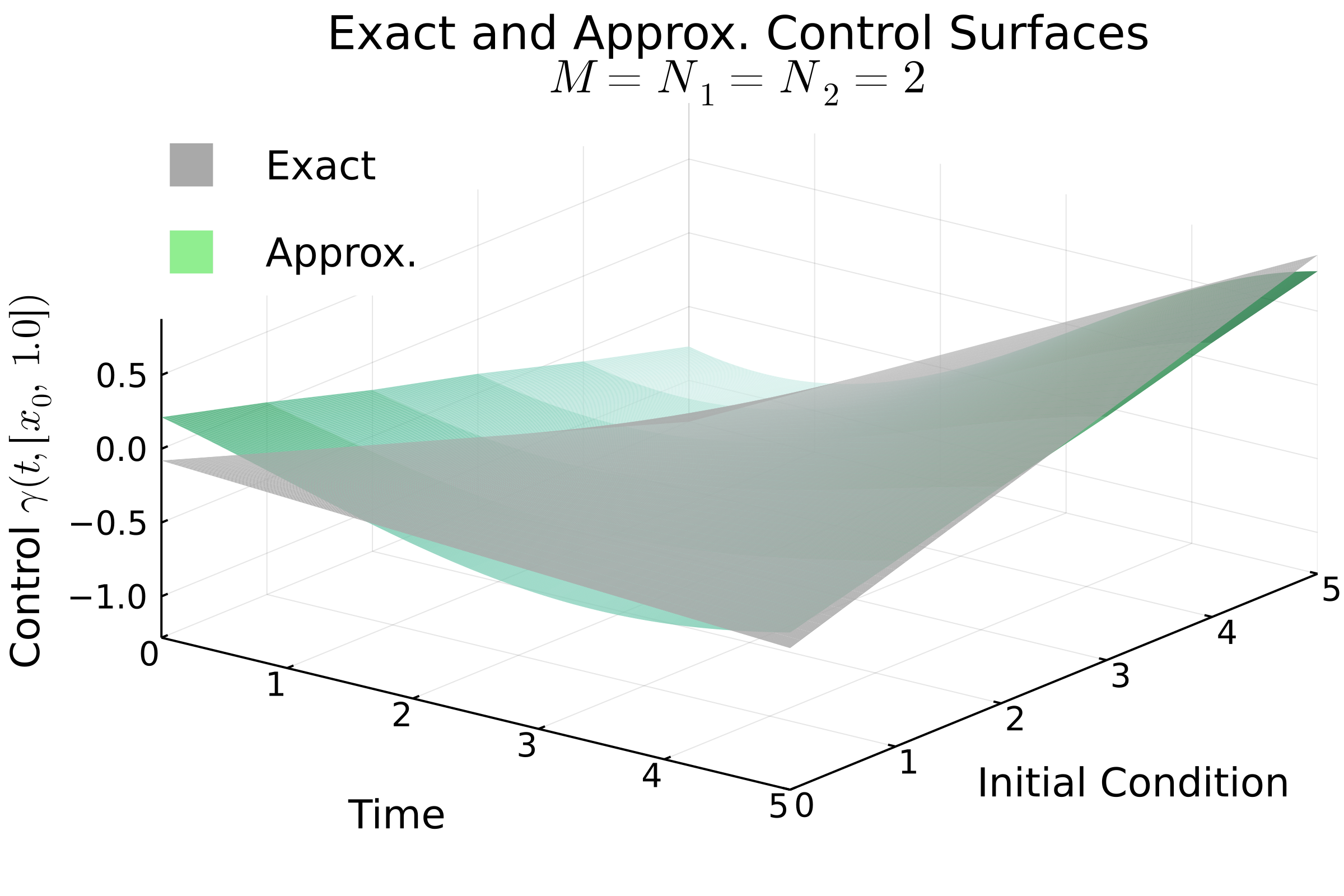}
    \qquad
    \includegraphics[width=0.6\textwidth]{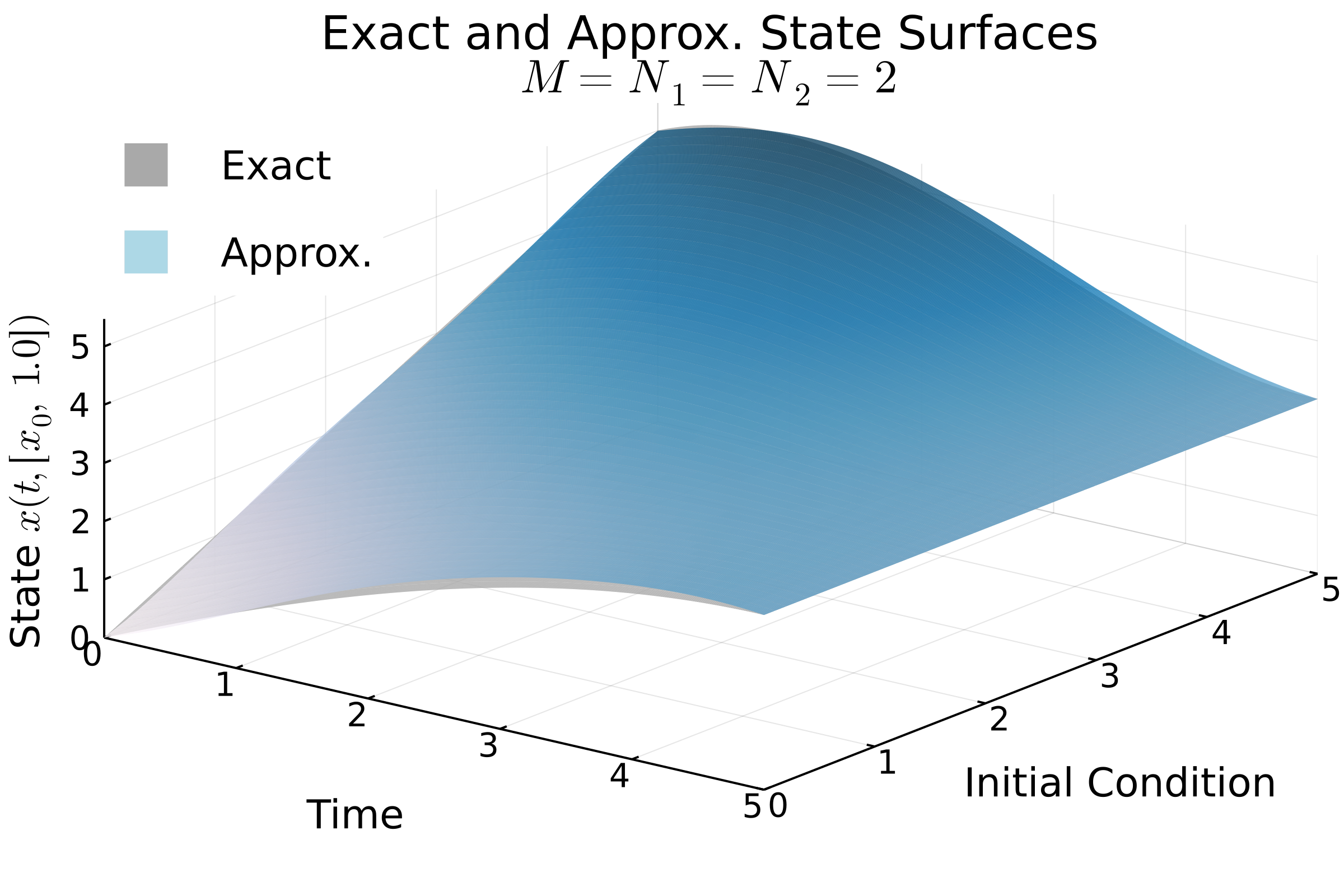}
    \caption{(Top) Slice at $\dot{x}=1.0$ of $4-$dimensional control surface approximation for $M=N_1=N_2=2$. (Bottom) Corresponding optimal state (position) surface.}
    \label{fig: toy2D v1}
\end{figure}

\begin{figure}[ht!]
\centering
    \centering
    \includegraphics[width=0.6\textwidth]{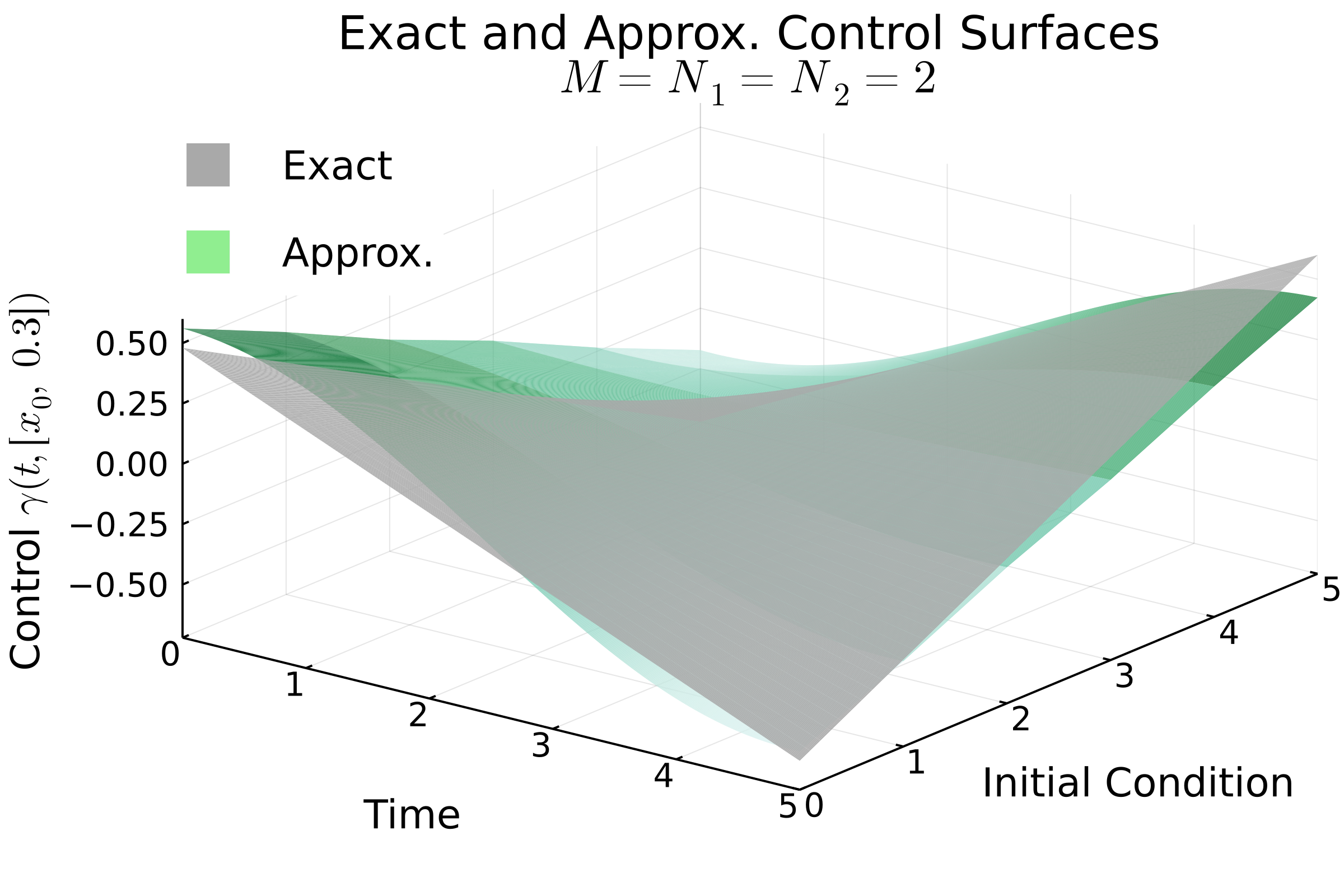}
    \qquad
    \includegraphics[width=0.6\textwidth]{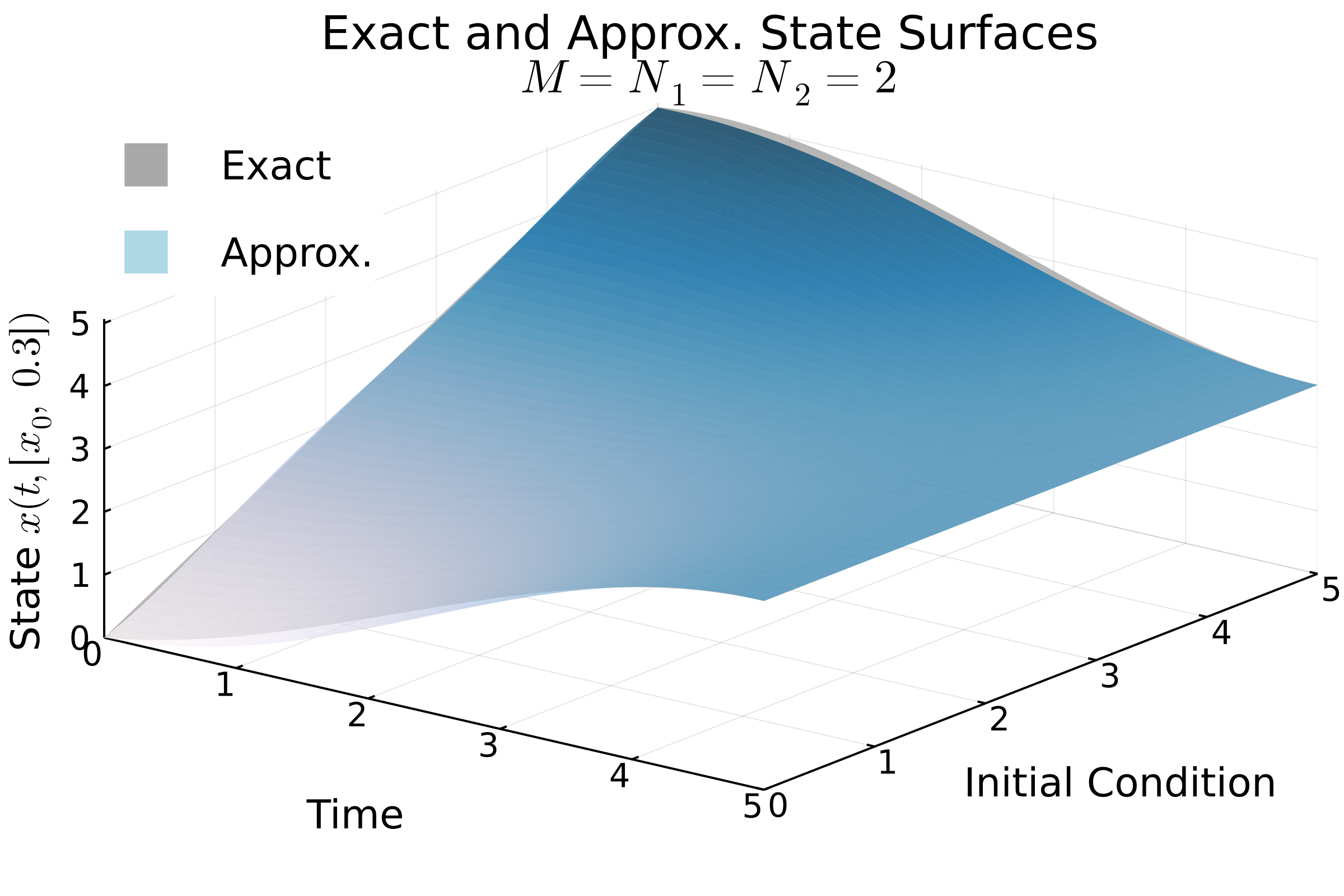}
    \caption{(Top) Slice at $\dot{x}=0.3$ of $4-$dimensional control surface approximation for $M=N_1=N_2=2$. (Bottom) Corresponding optimal state (position) surface.}
    \label{fig: toy2D v03}
\end{figure}

\begin{figure}[ht!]
\centering
    \centering
    \includegraphics[width=0.6\textwidth]{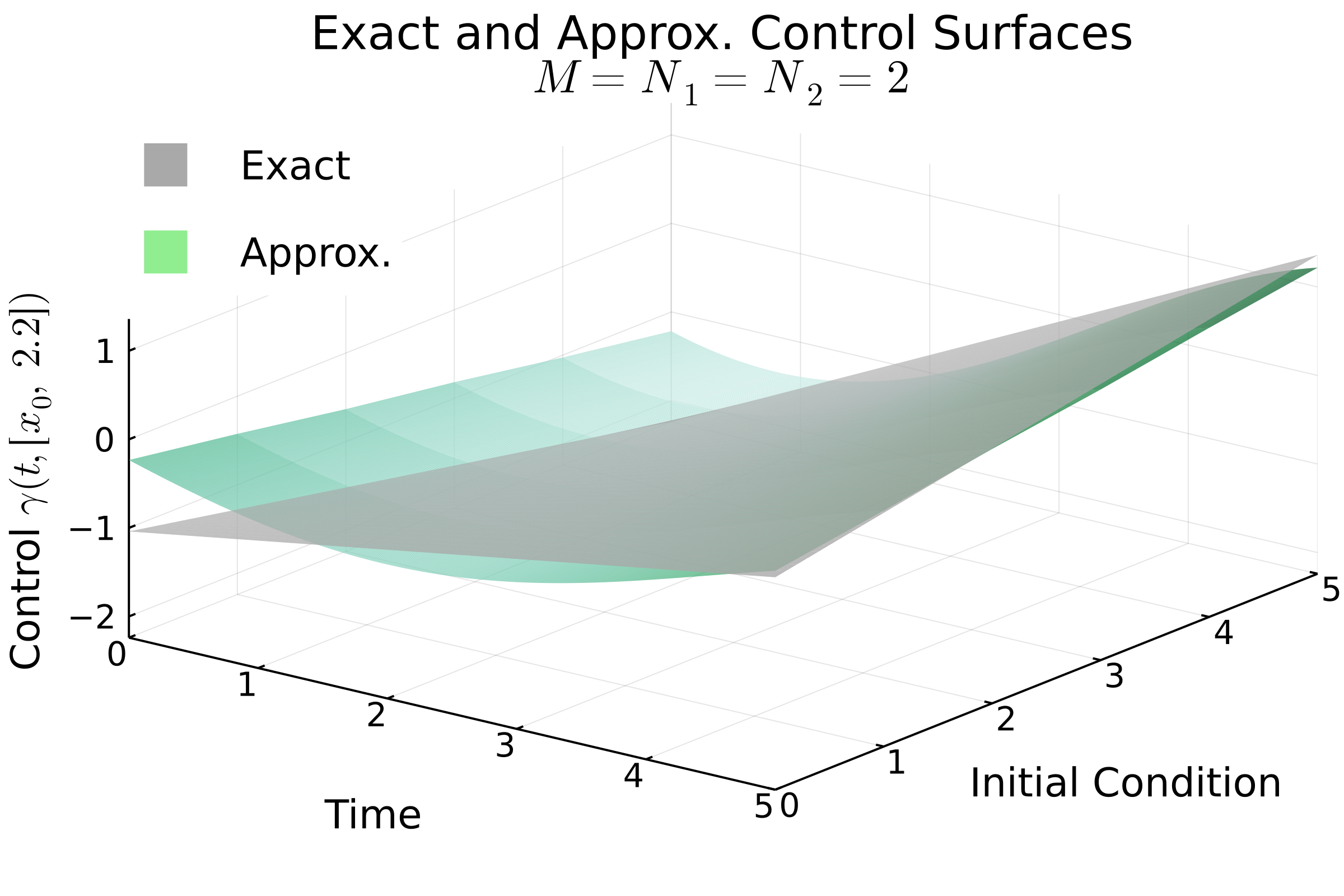}
    \qquad
    \includegraphics[width=0.6\textwidth]{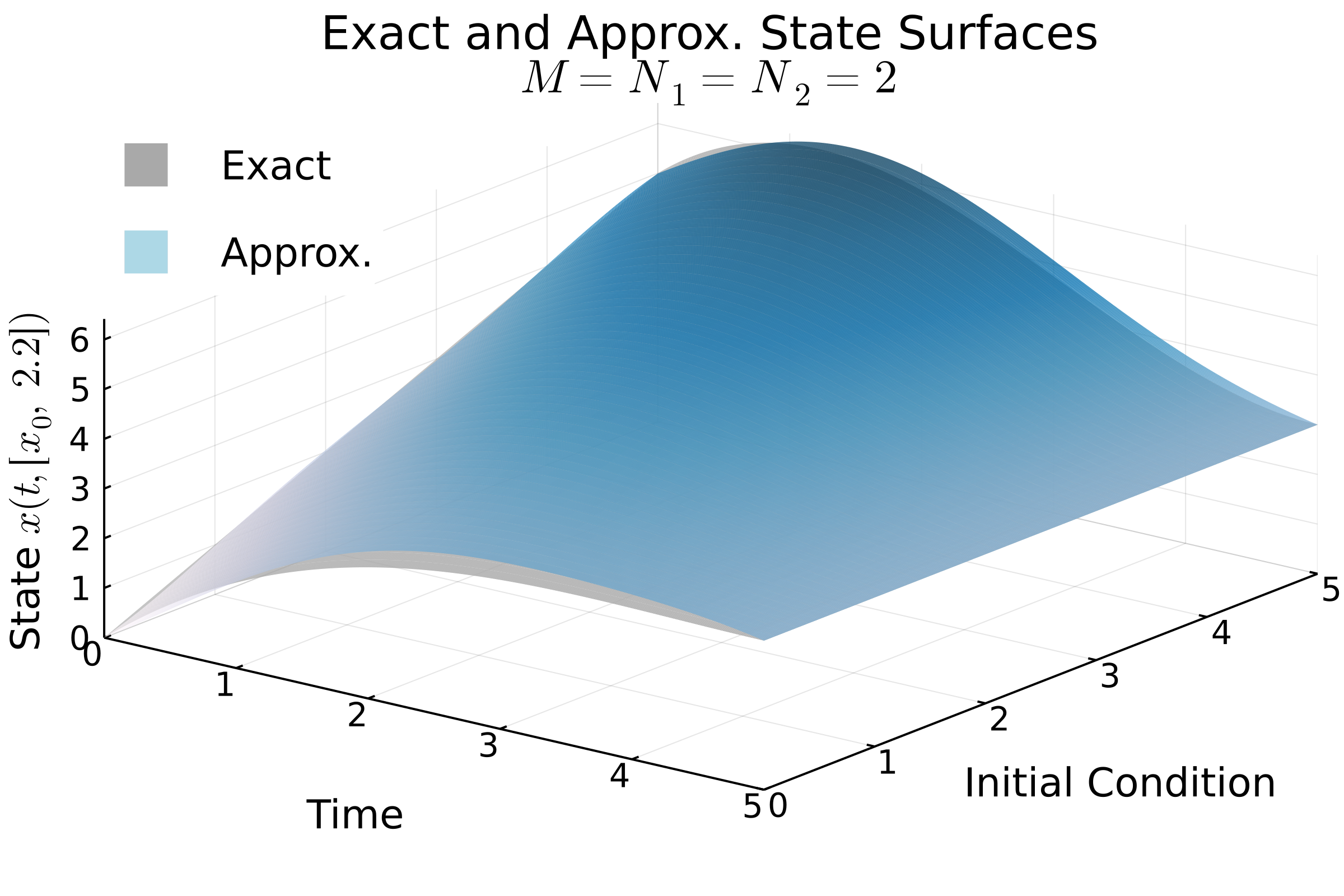}
    \caption{(Top) Slice at $\dot{x}=2.2$ of $4-$dimensional control surface approximation for $M=N_1=N_2=2$. (Bottom) Corresponding optimal state (position) surface.}
    \label{fig: toy2D v22}
\end{figure}

\cref{fig: toy2D v1,fig: toy2D v03,fig: toy2D v22} show the optimal control and trajectory approximations for three different values of velocity $\dot{x}$, respectively. We note that $\dot{x}=1$ belongs to the tuples in $X_0$, while the other two values, $\dot{x} = 0.3$ and $\dot{x} = 2.2$, do not, showing the generalization capability of the approximations with respect to initial state conditions tuples not in $X_0$.  

\subsection{Problem 2: Rock-Paper-Scissors} 
A particular kind of problem is that of OCP~(\ref{eqn: mainOCP}) for $d=3$, \textit{i.e.}, the famous rock-paper-scissors game. We recall that the functions $\mathbf{F}(\mathbf{u})$ and $\mathbf{G}(\mathbf{u})$ become,
\begin{equation}
    \begin{aligned}    
    &\mathbf{F}_i(\mathbf{u}) = u_i\left(\left(e_i - \mathbf{u}\right)^\top \mathbf{L}_3       \mathbf{u}\right) \\
    &\mathbf{G}_i(\mathbf{u}) = u_i\left(\left(e_i - \mathbf{u}\right)^\top \mathbf{M}_3       \mathbf{u}\right),    
    \end{aligned}
\end{equation}
where the circulant matrices $\mathbf{L}_3$ and $\mathbf{M}_3$ are,
\begin{equation}
        \mathbf{L}_3 = \begin{bmatrix} 0 & -1 & 1\\ 1 & 0 & -1\\ -1 & 1 & 0 \end{bmatrix}, \quad   \mathbf{M}_3 = \begin{bmatrix} 0 & 0  & 1\\ 1 & 0 &  0\\  0 & 1 & 0 \end{bmatrix}.
\end{equation}
Therefore, the parameterized OCP becomes,
\begin{equation}
    \begin{aligned}
        \min &\int_0^T \frac{1}{2}\lVert \mathbf{\hat{u}} - \mathbf{u}^* \rVert^2 + \frac{r}{2}\hat{\gamma}^2 \quad dt\\
        s.t. \;\; &\frac{d}{dt}\hat{\mathbf{u}} = \mathbf{F}(\hat{\mathbf{u}}) + \hat{\gamma}\mathbf{G}(\hat{\mathbf{u}})\\
        &\mathbf{u}(0) = \mathbf{u}_0 \in U_0, 
    \end{aligned}
    \label{eqn: mainOCP par}
\end{equation}
where $\hat{\gamma}$ and $\hat{\mathbf{u}}$ are Fourier approximations. Next, we consider the case with fixed initial condition $\mathbf{u}_0$, implying time-dependent approximations only, that is, Fourier approximations of the type $\hat{\gamma} = \hat{\gamma}(t, \Theta)$ and $\hat{\mathbf{u}} = \hat{\mathbf{u}}(t, \Theta_u)$. Second, the case with varying initial condition $\mathbf{u}_0$ in the unit simplex, 
\begin{equation}
\Delta_2 = \{ \mathbf{u} \in \mathbb{R}^3\;\vert\; \mathbf{1}^\top\mathbf{u} = 1,\; \mathbf{u} \ge \mathbf{0} \}
\label{eqn: 3D unit simplex}
\end{equation}
is considered. Note how this set defining the feasible state variables imposes another constraint to \cref{eqn: mainOCP par}, which is also included in the Augmented Lagrangian formulation, necessary to obtain the intended approximations. 

\subsubsection{Fixed Initial Condition} 
\begin{figure}[ht!]
\centering
    \centering
    \includegraphics[width=0.7\textwidth]{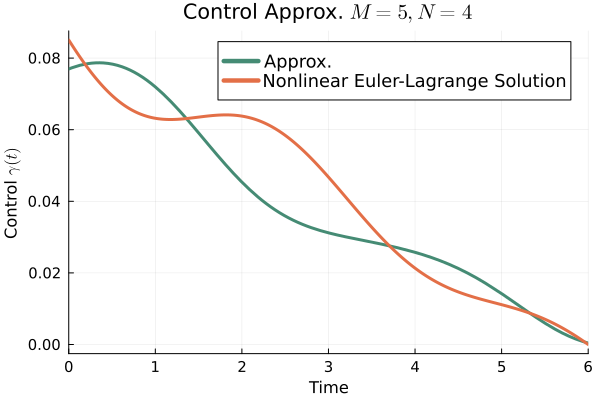}
    \qquad
    \includegraphics[width=0.7\textwidth]{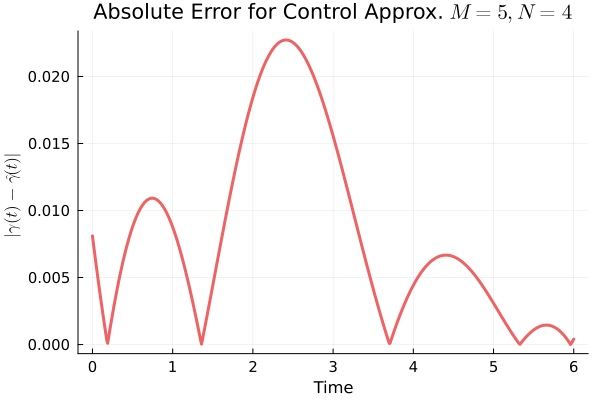}
    
    \caption{(Top) Optimal control approximation (green line) with $M=5$ and $N=4$ for fixed initial condition $\mathbf{u}_0 = [0.2333,\;0.3333,\;0.4333]$, and exact optimal control solution from the nonlinear Euler-Lagrange equation (orange). (Bottom) Corresponding absolute error.}
    \label{fig: Control Evol 1D M14}
\end{figure}

\begin{figure}[ht!]
\centering
    \centering
    \includegraphics[width=0.6\textwidth]{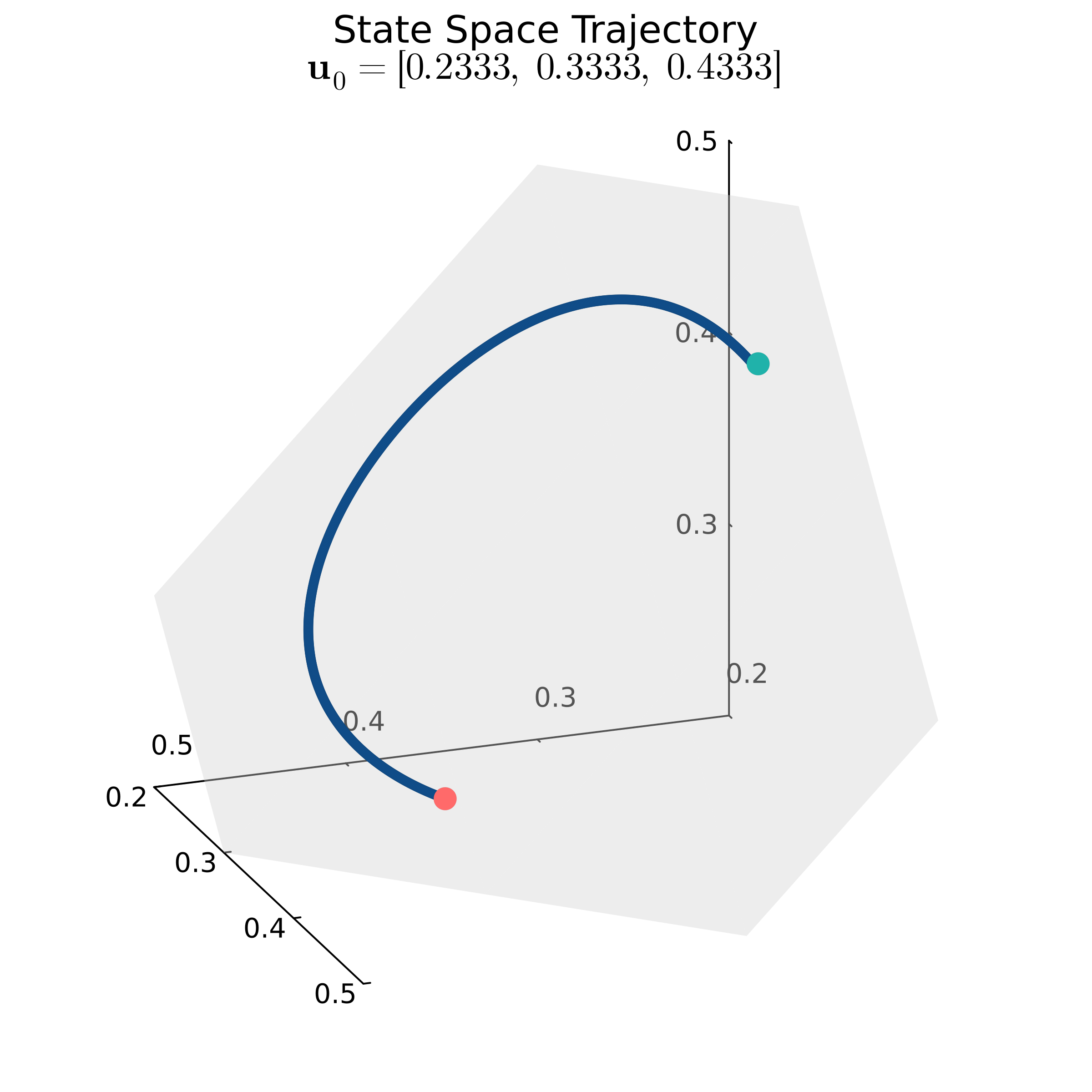}
    \qquad
    \includegraphics[width=0.6\textwidth]{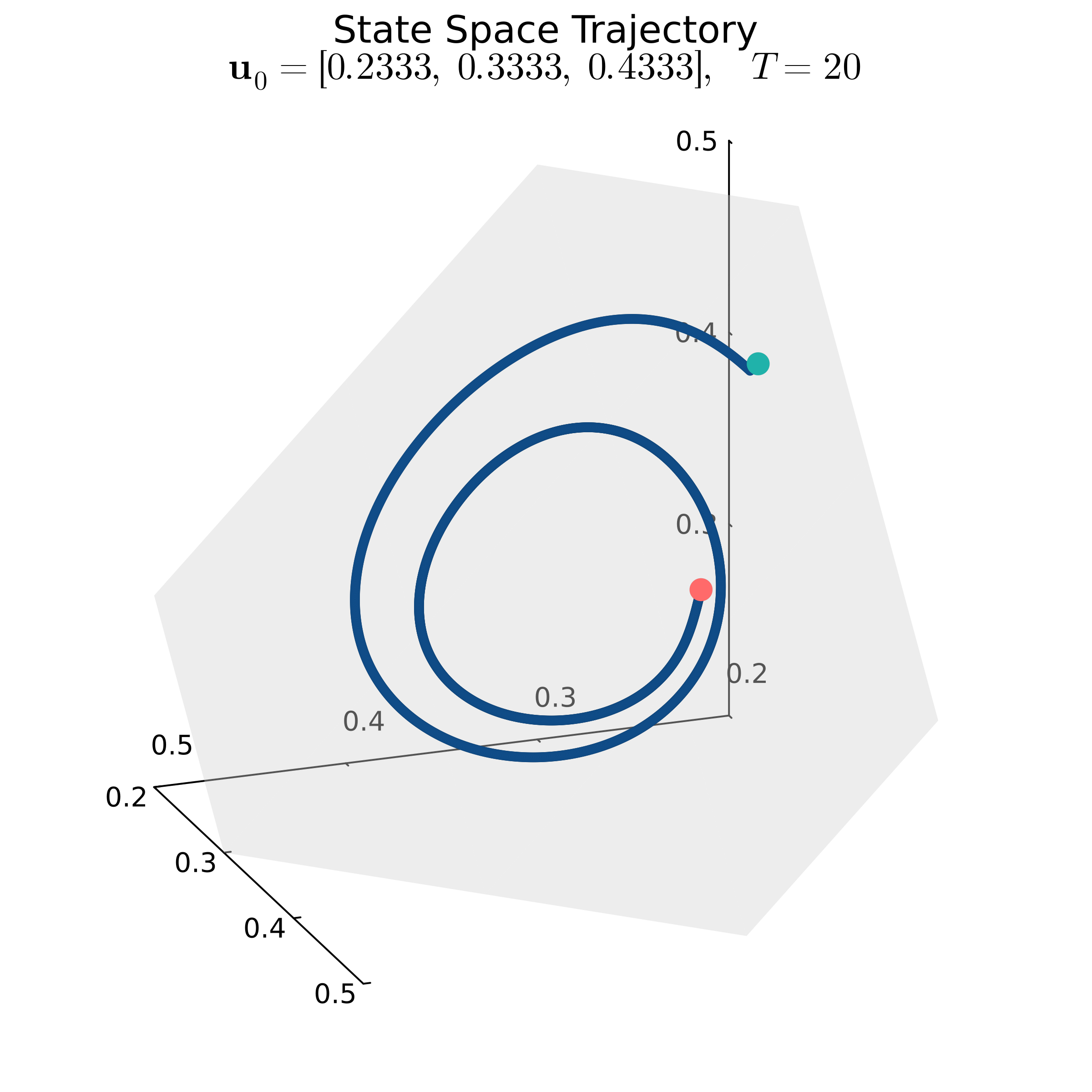}
    
    \caption{(Top) Corresponding optimal state trajectory (blue line) on the unit simplex (grey plane) following the approximating control of \cref{fig: Control Evol 1D M14}, starting at $\mathbf{u}_0$ (green) and ending at $\hat{\mathbf{u}}(T)$ (red). (Bottom) \textit{Idem} for when $T=20$, better depicting the spiralling trajectory toward $\mathbf{u}^*$.}
    \label{fig: Evol 1D M14}
\end{figure}

In OCP~(\ref{eqn: mainOCP par}), we set $U_0 = \{\begin{pmatrix} \frac{7}{30}, & \frac{1}{3}, & \frac{13}{30}\end{pmatrix}\}$, that is, the trajectory starting at a fixed initial condition in the $3-$dimensional unit simplex. This is a convenient choice of initial condition, allowing us to compare the results of this work to those presented in \cite{GF22}. A constraint enforcing the dynamics $\mathbf{F(\hat{u})} + \hat{\gamma}\mathbf{G(\hat{u})}$ is required, plus a constraint imposing \cref{eqn: 3D unit simplex} for each component in the state space. The Augmented Lagrangian in this case is formed by \crefrange{eqn:fpar}{eqn:h4par}. \cref{fig: Control Evol 1D M14} illustrates the approximated optimal control for a choice of $M=5$ and $N=4$ Fourier coefficients, alongside the exact optimal control derived from the solution of the Euler-Lagrange system of equations. \cref{fig: Evol 1D M14} shows the trajectory in state space following the approximated optimal control from \cref{fig: Control Evol 1D M14} for when $T=6$ and $T=20$. 

From \cref{tab: Evol1d}, we conclude that the number of outer-loops in \cref{alg: AL} seems to have a positive effect in decreasing $J^*_{\% error}$, as the smallest values for this measure were obtained for $\ell=100$ and $\ell=150$. In this regard, in all three experiments, adopting $\ell=100$ yielded approximately $J^*_{\% error} = 1.58\%$. These experiments also reveal that with respect to decreasing $J^*_{\% error}$, the number of coefficients $M$ and $N$ seems to have no significant effect. Additionally, we note how the case $\ell=150$ yields a worse $J^*_{\% error}$ when compared to the case in which $\ell = 100$, illustrating the fact that the choice of $\ell$ plays an important role in the convergence of \cref{alg: AL}, as discussed in \cite{bazaraa2013nonlinear, izmailov2007otimizaccao, NW99}. The impediment to attaining smaller values of $J^*_{\% error}$ might not be related to the choice of $\ell$ only. Indeed, this can have two other causes: first, not enough coefficients are being used in the truncated series. Second, \cref{alg: AL} might be converging to a local minimum, which is a common outcome in non-convex optimization. To address the first cause, more coefficients would have to be used in the truncated series, at the expense of more computational requirements.  \cref{sec: Error Bounds 1D,sec: Error Bounds 2D} will address the problem of bounding the mean square error for the optimal control approximation in terms of the number of coefficients used. To address the second potential cause, global optimization methods could be applied. It will also be necessary in future work to address the initialization of \cref{alg: AL} or even other choices of algorithm, such as stochastic gradient descent, known to better avoid convergence to local minima \cite{GB16}. 
\begin{table}[htbp]
    \centering
    \caption{Summary of computational results of Problem 2 for fixed $u(0)$. The constant parameters used in all the experiments are $k=2e+3$ and $\epsilon=10e-8$.}
    \label{tab: Evol1d}
    \begin{tabular}{ccccccc}
        \hline
        Method & $M = N + 1 $ & $\ell$ & $MSE_\gamma$ & $MAPE_\gamma$ & $MAE_\gamma$ & $J^{*}_{\%error}$ \\
        \hline
        \hline
        \multirow{7}{*}{LBFGS} & $4$ & $30$ & $1.10e-4$ & $502\%$ & $8.06e-3$ & $20.60\%$ \\
        & $4$  & $100$ & $1.40e-4$ &$64,882\%$& $1.00e-2$ & $1.58\%$ \\
        & $5$  & $30$ & $1.13e-4$ & $6,001\%$ & $8.12e-3$ & $20.60\%$ \\
        & $5$  & $100$ & $1.59e-4$ & $11,067\%$ & $1.05e-2$ & $1.58\%$ \\
        & $5$  & $150$ & $1.89e-4$ & $38,915\%$ & $1.13e-2$ & $1.72\%$\\
        & $6$  & $30$ & $1.12e-4$ & $1,802\%$ & $8.12e-3$ & $20.60\%$ \\
        & $6$  & $100$ & $1.57e-4$ & $2,162\%$ & $1.05e-2$ & $1.58\%$ \\
        & $7$  & $30$ & $1.12e-4$ & $3,345\%$ & $8.04e-3$ & $20.60\%$ \\
        \hline
    \end{tabular}
\end{table}
\clearpage
\section{Error Bounds for the Truncated Fourier Approximation}\label{sec: Error Bounds}
In this section, bounds for the mean square error (MSE) in approximating the optimal control (state) function with a truncated Fourier series are developed. The goal is to derive a MSE bound, which depends on the choice of the number of Fourier coefficients used in the truncated series. These results are fundamental when using the approximations developed in this work, as they provide us with the number of coefficients in the truncated series which guarantee the MSE of the approximation is not greater than a prescribed $\epsilon>0$. First, in \cref{sec: Error Bounds 1D}, we present a derivation for the fixed initial condition case. This derivation is inspired by \cite{G72}, which derives bounds for a periodic, single variable function over the interval $\left[-T/2, T/2\right]$. To the author's best knowledge, \cite{G72} offers the only bound derivation for a truncated Fourier approximation available in the literature. Lastly, in \cref{sec: Error Bounds 2D}, MSE bounds for the double Fourier approximation are provided, following a similar rationale to that in the single variable case. This latter result, to the best of the author's knowledge, is unprecedented in the literature. Nonetheless, a clear connection between the single and the two-variable case can be readily noted, suggesting a general MSE bound form for the $n-$dimensional approximation. 

\subsection{Error Bounds for the Fixed Initial Condition Case}\label{sec: Error Bounds 1D}
For the fixed initial condition case, we can provide error bounds for the truncated Fourier series approximation. To facilitate this development, we consider the Fourier series in complex exponential form as in \cref{eqn: complex-exp Fourier}. Error bounds for a periodic function on $[-\frac{T}{2}, \frac{T}{2}]$ are provided in \cite{G72}. Despite this being a rather intuitive result, here we show that the error bounds for a periodic function over $[0, 2T]$ are derived in the same manner as in the original case presented in \cite{G72}. We also provide a more detailed derivation than that presented in the original paper. Consider, for instance, the control function $\gamma(t)$ that we wish to approximate by $\hat{\gamma}_K(t)$, a truncated single variable Fourier series in complex exponential form, with $2K+1$ coefficients ($k=-K,-K+1,\dots,-1,0,1,\dots,K-1,K$). Consider the mean square error of such an approximation over $[0,2T]$, that is,
\begin{equation}
    MSE = \frac{1}{2T}\int_0^{2T} ( \gamma(t) - \hat{\gamma}_K(t) )^2 \; dt.
    \label{eqn: MSE}
\end{equation}
With $\hat{\gamma}_K(t) = \sum_{k=-K}^K c_k e^{\frac{i \pi k t}{T}}$, we can compute each coefficient as $c_k = \frac{1}{2T} \int_0^{2T} \gamma(t)e^{-\frac{i \pi k t}{T}} \; dt$ \cite{T76}, for $K\rightarrow\infty$. Furthermore, let us assume that the true optimal control is of bounded variation over $[0,2T]$, meaning that,
\begin{equation}
    V_{0}^{2T} [\gamma(t)] = \int_0^{2T} \vert \dot{\gamma}(t) \vert \; dt \le C.
    \label{eqn: bounded var}
\end{equation}
Before delving into the derivation of the bound, let us recall an important theorem from the theory of Fourier series. 
\begin{theorem}[Parseval's Theorem] For a periodic function $\gamma(t)$ over $[0,2T]$ with the Fourier series \cref{eqn: complex-exp Fourier}, the following equality holds 
\begin{equation}
\frac{1}{2T} \int_0^{2T} \gamma(t)^2 \; dt = \sum_{k= -\infty}^\infty \lvert c_k \rvert^2.      
\label{eqn: Parseval's}
\end{equation}
As a consequence, the MSE of the truncated Fourier series $\hat{\gamma}_K(t)$, subtracted from $\gamma(t)$, can be written as,
\begin{equation}
    MSE = \int_0^{2T} (\gamma(t) - \hat{\gamma}_K(t))^2 \; dt = 4T\sum_{k=K+1}^\infty \lvert c_k \rvert^2. \label{eqn: MSE Parseval}
\end{equation}
\label{thr: Parseval's}
\end{theorem}
\cref{eqn: MSE Parseval} provides a way to compute $MSE$ in \cref{eqn: MSE}. From our knowledge of Fourier series, we know that each $|c_k|$ term in \cref{eqn: MSE Parseval} can be computed as,
\begin{equation}
    \lvert c_k\rvert = \frac{1}{2T} \left\lvert \int_0^{2T} \gamma(t) e^{\frac{-i \pi k t }{T}}\;dt \right\rvert.
    \label{eqn: coeffi complex form}
\end{equation}
By noting that $\frac{d}{dt} e^{-i \pi k t /T} = (-i \pi k /T) e^{-i \pi k t /T}$, \cref{eqn: coeffi complex form} can be rewritten as a Riemann-Stieltjes integral,
\begin{equation}
    \lvert c_k \rvert = \frac{1}{\lvert k\rvert 2\pi} \left\lvert\int_0^{2T} \gamma(t) de^{\frac{-i\pi kt}{T}}\right\rvert.
\end{equation}
The integral in the absolute value above is now solved through integration by parts. Recalling that $\int u \; dv = uv - \int v \; du$, we have,
\begin{equation}
    \int_0^{2T} \gamma(t) de^{\frac{-i\pi kt}{T}} = \left. \gamma(t) e^{\frac{-i\pi kt}{T}}\right]_0^{2T} - \int_0^{2T} e^{\frac{-i\pi k t}{T}}d\gamma(t),  
\end{equation}
and,
\begin{equation}
    \left. \gamma(t) e^{\frac{-i\pi kt}{T}}\right\vert_0^{2T} = \gamma(2T)e^{-2i\pi k} - \gamma(0) = \gamma(2T) - \gamma(0) = 0.
\end{equation}
The second equality above follows from Euler's equation, which tells us that $e^{-i \pi} = -1$, and from the periodicity of $\gamma$ over $[0,2T]$. Thus, $\Pi_1^k e^{-i \pi} = e^{-i \pi k} = 1$ for even $k$ and $=-1$ for odd $k$. We then conclude that,
\begin{equation}
    \int_0^{2T} \gamma(t) de^{\frac{-i\pi kt}{T}} = - \int_0^{2T} e^{\frac{-i\pi k t}{T}}d\gamma(t).\label{eqn: equivalence}
\end{equation}
Now, we lay out some inequalities in order to find an upper bound for $|c_k|$. Substituting $d\gamma(t)$ in the above by $\dot{\gamma}(t)dt$, we have,
\begin{multline}
        \left\lvert \int_0^{2T} e^{\frac{-i \pi k t}{T}} \dot{\gamma}(t) \;dt  \right\rvert \le \int_0^{2T}  \left\lvert e^{\frac{-i \pi k t}{T}} \dot{\gamma}(t) \right\rvert \;dt
         = \\
         \int_0^{2T} \left\lvert e^{\frac{-i \pi k t}{T}} \right\rvert \left\lvert \dot{\gamma}(t)\right\rvert\; dt 
         = \int_0^{2T} \left\lvert \dot{\gamma}(t) \right\rvert \; dt 
      \le C,
        \label{eqn: inequalities}
\end{multline}
where the second equality is justified since $\left\lvert e^{\frac{-i \pi k t}{T}}\right\rvert = 1$, and the last inequality is simply \cref{eqn: bounded var}. The result in \cref{eqn: inequalities} implies,
\begin{equation}
    \lvert c_k \rvert = \frac{1}{\lvert k \rvert 2\pi} \left\lvert \int_0^{2T} e^{\frac{-i \pi \lvert k\rvert t}{T}} \dot{\gamma}(t) \;dt  \right\rvert \le \frac{1}{\lvert k\rvert 2\pi} V_0^{2T} [\gamma(t)] \le \frac{C}{\lvert k \rvert 2\pi}.
\end{equation}
Invoking \cref{eqn: MSE Parseval}, we have,
\begin{equation}
    MSE \le 4T\sum_{k=K+1}^\infty \frac{C^2}{4k^2\pi^2} = \frac{TC^2}{\pi^2} \sum_{k=K+1}^\infty \frac{1}{k^2}.
    \label{eqn: MSE inequality}
\end{equation}
To compute the summation term in the expression above, we note that the following inequality holds,
\begin{equation}
    \sum_{k=K+1}^\infty \frac{1}{k^2} \le \int_K^\infty \frac{1}{x^2} \; dx = \left. -\frac{1}{x} \right]_K^\infty = \frac{1}{K}.
\end{equation}
Thus, we conclude that,
\begin{equation}
    MSE \le \frac{T{C}^2}{\pi^2K}.
    \label{eqn: MSE bound}
\end{equation}
This result provides us with a number of coefficients $K$ used in the truncated Fourier approximation of $\gamma(t)$, under the conditions already discussed. For example, if we wish to obtain an approximation with a $MSE \le \epsilon$, we set $\epsilon = \frac{T{C}^2}{\pi^2 K}$ and solve for $K$. Any approximation involving at least $2K+1$ coefficients is guaranteed to achieve a $MSE$ less than the specified value of $\epsilon$. It is important to remark that this result does not eliminate the possibility of obtaining an approximation that respects $MSE \le \epsilon$ when the number of coefficients is not given by \cref{eqn: MSE bound}. We also note the case in which $k=0$ is not relevant given the summation limit in \cref{eqn: MSE inequality}. 

\subsection{Error Bounds for the Varying Initial Condition or Two-variable Fourier Series Case}\label{sec: Error Bounds 2D}
The same analysis performed in the previous paragraphs can be extended for the case involving a Fourier approximation of two variables, these being time and initial condition. For this, we will consider the control approximation $\hat{\gamma}_{KL}(t, u_0)$ with $\gamma(t,u_0)$ denoting both the actual optimal control and its full Fourier representation for simplified notation. We also assume that $\gamma$ has a periodic extension over the domain defined by $[0,2T]\times \mathcal{U}_0$, where $\mathcal{U}_0 = [0,2U]$ for the purposes of this analysis. Again, the mean square error obtained by the truncated approximation is,
\begin{equation}
    MSE = \frac{1}{4TU}\int_0^{2U} \int_0^{2T} \left(\gamma(t, u_0) - \hat{\gamma}_{KL}(t, u_0) \right)^2 \; dt\;du_0.
\end{equation}
The fundamental assumption on $\gamma(t,u_0)$ is that it is a function of bounded total variation, that is,
\begin{equation}
    V_{0,0}^{2T,2U}\left[\gamma(t,u_0)\right] = \int_0^{2T} \int_0^{2U} \left\lvert \nabla\gamma(t,u_0)\right\rvert \;du_0\;dt \le C. 
\end{equation}
From now on, we proceed as before. The truncated Fourier approximation is,
\begin{equation}
        \hat{\gamma}_{KL}(t, u_0) = \sum_{k=-K}^K\sum_{l=-L}^L c_{k,l} e^{\frac{ik\pi t}{T}} e^{\frac{il\pi u_0}{U}},
\end{equation}
and the closed-form coefficients of the full (complex form) Fourier series of $\gamma(t,u_0)$ are given by,
\begin{equation}
    c_{k,l} = \frac{1}{4TU} \int_0^{2U} \int_0^{2T} \gamma(t,u_0) e^{\frac{-ik\pi t}{T}} e^{\frac{-il\pi u_0}{U}} \; dt\;du_0. 
    \label{eqn: 2D Fourier coeffis}
\end{equation}
Invoking Parseval's theorem, the mean square error can be written as,
\begin{equation}
    MSE = 16TU \sum_{k=K+1}^\infty\sum_{l=L+1}^\infty \lvert c_{k,l} \rvert^2. 
\end{equation}
In computing $\lvert c_{k,l}\rvert$, we write,
\begin{equation}
    \lvert c_{k,l}\rvert  = \frac{1}{4TU} \left\lvert \int_0^{2U} \int_0^{2T} \gamma(t,u_0) e^{\frac{-ik\pi t}{T}} e^{\frac{-il\pi u_0}{U}} \; dt\;du_0 \right\rvert .
\end{equation}
It is convenient in this case to not use the Riemann-Stieltjes form. Using the same argument as in \cref{eqn: equivalence}, the inner integrand of \cref{eqn: 2D Fourier coeffis} becomes,
\begin{multline}
    \int_0^{2T} \gamma(t, u_0)e^{\frac{-ik\pi t}{T}}\;dt
    = \\
    \frac{T}{ik\pi}\left[ \left[\gamma(0, u_0) - \gamma(2T, u_0)\right] + \int_0^{2T} \frac{\partial \gamma(t, u_0)}{\partial t}e^{\frac{-ik\pi t}{T}}\;dt \right]= I_T.
\end{multline}
We call the result obtained above of $I_T$, and, in computing the double integral term in $c_{k,l}$, we are now left with $\int_0^{2U} I_T \; e^{\frac{-il\pi u_0}{U}}\;du_0$. Breaking down this integral into three components, we obtain,
\begin{multline}
    \int_0^{2U} \gamma(0,u_0)e^{\frac{-il\pi u_0}{U}}\;du_0 = \\
    \frac{U}{il \pi}\left[\left[\gamma(0,0) - \gamma(0,2U)\right] + \int_0^{2U} \frac{\partial \gamma(t,u_0)}{\partial u_0}e^{\frac{-il\pi u_0}{U}}\;du_0\right]
\end{multline}
\begin{multline}
    \int_0^{2U} \gamma(2T,u_0)e^{\frac{-il\pi u_0}{U}}\;du_0 =\\
    \frac{U}{il \pi}\left[\left[\gamma(2T,0) - \gamma(2T,2U)\right] + \int_0^{2U} \frac{\partial \gamma(t,u_0)}{\partial u_0}e^{\frac{-il\pi u_0}{U}}\;du_0\right].
\end{multline}
From the periodicity of $\gamma$, we know that,
\begin{equation}
\gamma(0,0) - \gamma(0,2U) = \gamma(2T,0) - \gamma(2T,2U) = 0.
\end{equation}
And, finally, the third component is bounded above as,
\begin{equation}
\begin{aligned}
    \int_0^{2U} \int_0^{2T} \frac{\partial \gamma(t, u_0)}{\partial t}e^{\frac{-il\pi u_0}{U}}e^{\frac{-ik\pi t}{T}}\;dt\;du_0 &\le \int_0^{2U}\int_0^{2T} \left\lvert \frac{\partial \gamma(t,u_0)}{\partial t} e^{\frac{-il\pi u_0}{U}}e^{\frac{-ik\pi t}{T}} \right\rvert\;dt\;du_0 \\
    &\le \int_0^{2U}\int_0^{2T} \left\lvert \frac{\partial \gamma(t,u_0)}{\partial t}\right\rvert \left\lvert e^{\frac{-il\pi u_0}{U}}\right\rvert \left\lvert e^{\frac{-ik\pi t}{T}}\right\rvert \;dt\;du_0 \\
    &\le \int_0^{2U}\int_0^{2T} \left\lvert \nabla\gamma(t,u_0)\right\rvert \;dt\;du_0 \\
    &\le C.
\end{aligned}
\label{eqn: partIII}
\end{equation}
We conclude by writing,
\begin{equation}
\begin{aligned}
\lvert c_{k,l} \rvert \le &\frac{1}{4TU}\left\lvert \frac{T}{ik\pi}\frac{U}{il\pi}2C \right\rvert
= \frac{C}{\lvert k\rvert\lvert l\rvert2\pi^2}.
\end{aligned}
\end{equation}
Now we have obtained an upper bound for $\lvert c_{k,l} \rvert$. Using Parseval's theorem again, we have,
\begin{equation}
    MSE \le 16TU \sum_{k=K+1}^\infty\sum_{l=L+1}^\infty \frac{C^2}{k^2 l^2 4\pi^4} = \frac{4TUC^2}{\pi^4} \sum_{k=K+1}^\infty \frac{1}{k^2} \sum_{l=L+1}^\infty \frac{1}{l^2} \le \frac{4TUC^2}{\pi^4KL}.
\end{equation}
As we did before for the single variable case, if we set a $\epsilon > 0$ value which the $MSE$ of the truncated approximation is not supposed to exceed, we simply solve  $\epsilon = \frac{4TUC^2}{\pi^4KL}$ for $KL$.  Thus, for any values of $K$ and $L$ such that $KL$ is equal or greater than the value provided by the equation will guarantee that $MSE \le \epsilon$. 

\section{Conclusions}\label{sec: Conclusions}
In this work, we proposed a multidimensional truncated Fourier approximation to optimal control and trajectory hypersurfaces. We developed an Augmented Lagrangian algorithm that translates the OCP into an unconstrained optimization problem whose gradient can be computed via automatic differentiation, making our solution approach appealing from the machine learning perspective. Two problems were used to test our approach: Problem 1 considered a simple quadratic OCP describing the movement of a particle under Newton's law of motion. Problem 2 was a more elaborated OCP arising from the replicator dynamics in an odd-circulant game. Our computational experiments showed that, for the two-dimensional case of Problem 1 (time and initial position with fixed initial velocity), percentage errors for the optimal cost less than $2.5\%$ were obtained with relatively few coefficients in the truncated series. For the three-dimensional case, the best percentage error obtained was of approximately $20\%$. For Problem 2, only the one dimensional case (fixed initial condition) was considered. The results reported showed that a percentage error of less than $1.6\%$ was possible with a few coefficients in the truncated series. For higher dimensional cases, the computational requirements showed to be prohibitive, and this will be the matter of future investigation. Initialization of the proposed algorithm needs to be addressed, preventing the solution from getting stuck at a local minimum. Global optimization methods will also be considered in an attempt to circumvent local minima. Bounds for the MSE for the one and two-dimensional cases were also derived, this latter result being unprecedented in the literature. Nonetheless, our approach proved to be effective in translating the OCPs into an unconstrained optimization problem that can be solved using modern machine learning automatic differentiation tools and out-of-shelf unconstrained optimization solvers. The classical results of convergence in Fourier analysis together with the bounds derived in this work suggest that our method is successful in approximating optimal control and trajectories for a relatively broader class of OCPs.   

\clearpage

\printbibliography

\end{document}